\documentclass[10pt,leqno]{article}

\usepackage{amsmath,amssymb,amsthm,mathrsfs,dsfont}

\usepackage[margin=3cm]{geometry} 

\usepackage{titlesec,hyperref}

\usepackage{color}

\usepackage{fancyhdr}
\titleformat{\subsection}{\it}{\thesubsection.\enspace}{1pt}{}

\newtheorem{theo}{Theorem}[section]
\newtheorem{lemm}[theo]{Lemma}
\newtheorem{defi}[theo]{Definition}
\newtheorem{coro}[theo]{Corollary}
\newtheorem{prop}[theo]{Proposition}
\newtheorem{rema}[theo]{Remark}
\numberwithin{equation}{section}

\allowdisplaybreaks 

\begin{document}
\title{Well-posedness for the FENE dumbbell model of polymetric flows in Besov spaces
\hspace{-4mm}
}

\author{Wei Luo$^1$
\quad Zhaoyang Yin$^2$ \\[10pt]
Department of Mathematics, Sun Yat-sen University,\\
510275, Guangzhou, P. R. China.\\[5pt]
}
\footnotetext[1]{Email: \it luowei23@mail2.sysu.edu.cn}
\footnotetext[2]{Email: \it mcsyzy@mail.sysu.com.cn}
\date{}
\maketitle
\hrule

\begin{abstract}
In this paper we mainly investigate the Cauchy problem of the finite extensible nonlinear elastic (FENE) dumbbell model with dimension $d\geq2$. We first proved the local well-posedness for the FENE model in Besov spaces by using the Littlewood-Paley theory. Then by an accurate estimate we get a blow-up criterion. Moreover, if the initial data is perturbation around equilibrium, we obtain a global existence result. Our obtained results generalize recent results in \cite{Masmoudi.W}.\\

\vspace*{5pt}

\vspace*{5pt}
\end{abstract}

\vspace*{10pt}

\tableofcontents

\section{Introduction}
   In this paper we consider the finite extensible nonlinear elastic (FENE) dumbbell model \cite{Bird}:
   \begin{align}
\left\{
\begin{array}{ll}
\partial_{t}u+(u\cdot\nabla)u-\nu\Delta{u}+\nabla{P}=div\tau, ~~~~~~~div u=0,\\[1ex]
\partial_{t}\psi+(u\cdot\nabla)\psi=div_{R}[-\nabla{u}\cdot{R}\psi+\beta\nabla_{R}\psi+\nabla\mathcal{U}\psi],  \\[1ex]
\tau_{ij}=\int_{B}R_{i}\otimes\nabla_{j}\mathcal{U}\psi dR, \\[1ex]
u|_{t=0}=u_{0}, \psi|_{t=0}=\psi_{0} ,\\[1ex]
(\beta\nabla_{R}\psi+\nabla\mathcal{U}\psi)\cdot{n}=0 ~~~~ on ~~~~ \partial B(0,R_{0}) .\\[1ex]
\end{array}
\right.
\end{align}
In (1.1)~~$\psi(t,x,R)$ denotes the distribution function for the internal configuration and $u(t,x)$ stands for the velocity of the polymeric liquid, where $x\in\mathbb{R}^{d}$ and $d\geq2$ means the dimension. Here the polymer elongation $R$ is bounded in ball $ B=B(0,R_{0})$ of $\mathbb{R}^{d}$ which means that the extensibility of the polymers is finite. $\beta$ is a constant related to the temperature and $\nu>0$ is the viscosity of the fluid. $\tau$ is an additional stress tensor and $P$ is the pressure.

   This model describes the system coupling fluids and polymers. The system is of great interest in many branches of physics, chemistry, and biology, see \cite{Bird,Masmoudi.W}. In this model, a polymer is idealized as an "elastic dumbbell" consisting of two "beads" joined by a spring that can be modeled by a vector $R$. At the level of liquid, the system couples the Navier-Stokes equation for the fluid velocity with a Fokker-Planck equation describing the evolution of the polymer density. This is a micro-macro model (For more details, one can refer to  $\cite{Bird}$, $\cite{Masmoudi.W}$ and $\cite{Masmoudi.G}$).

%
%
%


In the paper we will take $\beta=1$ and $R_{0}=1$.
Notice that $(u,\psi)$ with $u=0$ and $$\psi_{\infty}(R)=\frac{e^{-\mathcal{U}(R)}}{\int_{B}e^{-\mathcal{U}(R)}dR},$$
is a stationary solution of (1.1). Thus we can rewrite (1.1) for the following system:
\begin{align}
\left\{
\begin{array}{ll}
\partial_{t}u+(u\cdot\nabla)u-\nu\Delta u+\nabla{P}=div\tau,  ~~~~~~~div u=0,\\[1ex]
\partial_{t}\psi+(u\cdot\nabla)\psi=div_{R}[-\nabla{u}\cdot{R}\psi+\psi_{\infty}\nabla_{R}\frac{\psi}{\psi_{\infty}}],  \\[1ex]
\tau_{ij}=\int_{B}R_{i}\otimes\nabla_{j}\mathcal{U}\psi dR, \\[1ex]
u|_{t=0}=u_{0}, \psi|_{t=0}=\psi_{0} ,\\[1ex]
\psi_{\infty}\nabla_{R}\frac{\psi}{\psi_{\infty}}\cdot{n}=0 ~~~~ on ~~~~ \partial B(0,R_{0}) .\\[1ex]
\end{array}
\right.
\end{align}
   Moreover the potential $\mathcal{U}(R)=-k\log(1-|R|^{2})$ for some constant $k>0$. We have to add a boundary condition to insure the conservation of $\psi$, namely,  $(-\nabla{u}\cdot{R}\psi+\psi_{\infty}\nabla_{R}\frac{\psi}{\psi_{\infty}})\cdot n=0$.
    The second equation in (1.2) can be understood in the weak sense: for any function $g(R)\in C^{1}(B)$, we have
   $$\partial_{t}\int_{B}g\psi dR+(u\cdot\nabla)\int_{B}g\psi dR
   =-\int_{B}\nabla_{R}g[-\nabla{u}\cdot{R}\psi+\psi_{\infty}\nabla_{R}\frac{\psi}{\psi_{\infty}}]dR.$$
\begin{defi}
Assume that $u_{0}\in S'(\mathbb{R}^{d})$ and $\psi_{0}\in S'(\mathbb{R}^{d};D'(B))$. A couple of functions $(u,\psi)\in C([0,T];S'(\mathbb{R}^{d}))\times C([0,T];S'(\mathbb{R}^{d};D'(B)))$ with $div~u=0$ is called a solution for (1.2) if for each $(v,\phi)\in C^{1}([0,T];S(\mathbb{R}^{d}))\times C^{1}([0,T];S(\mathbb{R}^{d});C^{\infty}(B))$ with $v(T)=0,~\phi(T)=0$ we have
\begin{align}
&\int^{T}_{0}\int_{R^{d}}u\partial_{t}v+(u\otimes u):\nabla v-\nu u\Delta v+P\cdot div~v=\int^{T}_{0}\int_{R^{d}}\tau:v+\int_{R^{d}}u_{0}v_{0}, \\
&\int^{T}_{0}\int_{R^{d}\times B}\psi\partial_{t}\phi+u\psi\cdot\nabla_{x}\phi =\int^{T}_{0}\int_{R^{d}\times B}[-\nabla{u}\cdot{R}\psi+\psi_{\infty}\nabla_{R}\frac{\psi}{\psi_{\infty}}]\cdot\nabla_{R}\phi +\int_{R^{d}\times B}\psi_{0}\phi_{0}.
\end{align}
\end{defi}
   Let us mention that the earliest local well-posedness for (1.1) was established by Renardy in \cite{Renardy}, where the author considered the Dirichlet problem with $d=3$ for smooth boundary and proved local existence for (1.1) in $\bigcap\limits^4_{i=0} C^i([0,T);H^{4-i})\times \bigcap\limits^3_{i=0}\bigcap\limits^{3-i}_{j=0} C^i([0,T);H^j)$ with potential $\mathcal{U}(R)=(1-|R|^2)^{1-\sigma}$ for $\sigma>1$. Later, Jourdain, Leli$\grave{e}$vre, and
Le Bris \cite{Jourdain} proved local existence of a stochastic
differential equation with potential $\mathcal{U}(R)=-k\log(1-|R|^{2})$ in the case $k>3$ for a Couette flow. Zhang and Zhang \cite{Zhang-P} proved local well-posedness of (1.1) with $d=3$ in $\bigcap\limits^2_{i=0} H^i([0,T);H^{4-2i})\times \bigcap\limits^1_{i=0} H^i([0,T);H^{3-2i})$ for $k > 38$. Lin, Zhang, and Zhang \cite{F.Lin} proved global well-posedness of (1.2) with $d=2$ for $ k > 6 $ in $C([0,T];H^s)\times C([0,T];H^s(R^2;H^1_0(D)))$, where $s\geq3$. Masmoudi \cite{Masmoudi.W} proved local well-posedness of (1.2) in $C([0,T);H^s)\times C([0,T);H^s(R^d; L^2))$ and global well-posedness of (1.2) when the initial data is perturbation around equilibrium for $k>0$. In the case $d=2$, the author  \cite{Masmoudi.W} obtained a global result for $k>0$.  Kreml and Pokorn$\acute{y}$ \cite{Kreml} proved local well-posedness in Sobolev spaces $W^{1,p}$ with $p>d$. Recently, global existence of weak solutions in $L^2$ was proved by Masmoudi \cite{Masmoudi.G} under some entropy conditions.

To our best knowledge, there were no results about the well-poesdness of (1.1) in Besov spaces. In this paper we investigate the well-posedness of (1.1) in Besov spaces $B^s_{p,r}$, which requires more elaborate techniques. In FENE model (1.1), the most difficult term is the additional stress tensor, however, Mousmoudi \cite{Masmoudi.G} proved a lot of useful lemmas to deal with this term. By using Mousmoudi's lemmas, we can easily get a corollary to solve the problem. Thus the remain difficulties are the product term and the pressure term. In order to obtain the well-poesdness of (1.1) in Besov spaces, one can apply $\Delta_j$ to (1.2) and get a localization of the equations. The product term leads to the commutator. If $p\neq2$, the energy method doesn't work. However, by using the Littlewood-Paley theory and Bony's decomposition, we split the commutator into 8 terms and for each term one can easily deal with by the basic H\"{o}lder's inequality. Then we obtain the commutator estimates which lead us to obtain a priori estimates. In order to deal with the pressure, one can apply $div$ to (1.2), since $div~u=0$ it follows that $p$ satisfies an elliptic equation. There is an explicit formula giving the pressure in terms of the velocity field, then by using some techniques in Fourier analysis, we can deal with the pressure term. Thanks to the viscosity coefficient $\nu>0$, the energy decays as time grows. If the $H^s$-norm of initial data is small, one can get the $H^s$-norm of the corresponding solution to (1.1) is smaller than that of initial data for any $t>0$. Then by an iteration argument, one can obtain the global well-poesdness of (1.1) in $H^s$. But for Bosov spaces $B^s_{p,r}$, if $p\neq2$, we can't obtain this property. However, we can use a continuous argument mentioned in \cite{B.C} to show that if the initial data is small then the corresponding solution is uniformly bounded by the initial data independent with $t$, which leads to the global well-poesdness of (1.1) in $B^s_{p,r}$.\\

  The paper is organized as follow. In Section 2 we introduce some notations and our main results. In Section 3 we give some preliminaries which will be used in this paper. In Section 4 we investigate the linear problem of (1.1) and give some a priori estimates for solutions to (1.1). In Section 5 we prove the local well-posedness of (1.1) by using approximate argument. Section 6 is devoted to the study of a blow-up criterion. In Section 7 we prove the global well-posedness of (1.1) by a contradiction argument.

\section{Notations and main results}
  In this section we introduce our main results and the notations that we shall use throughout the paper.

For $p\geq1$, we denote by $\mathcal{L}^{p}$ the space
$$\mathcal{L}^{p}=\big\{\psi \big|\|\psi\|^{p}_{\mathcal{L}^{p}}=\int \psi_{\infty}|\frac{\psi}{\psi_{\infty}}|^{p}dR<\infty\big\}.$$

  We will use the notation $L^{p}_{x}(\mathcal{L}^{q})$ to denote $L^{p}[\mathbb{R}^{d};\mathcal{L}^{q}]:$
$$L^{p}_{x}(\mathcal{L}^{q})=\big\{\psi \big|\|\psi\|_{L^{p}_{x}(\mathcal{L}^{q})}=(\int_{\mathbb{R}^{d}}(\int_{B} \psi_{\infty}|\frac{\psi}{\psi_{\infty}}|^{q}dR)^{\frac{p}{q}}dx)^{\frac{1}{p}}<\infty\big\}.$$
  Next we introduce the Littlewood-Paley decomposition and Besove spaces (see \cite{B.C.D} for more details).

Let $\mathcal{C}$ be the annulus $\{\xi\in\mathbb{R}^{d}\big|\frac{3}{4}\leq|\xi|\leq\frac{8}{3}\}.$ There exist radial functions $\chi$ and $\varphi$, valued in the interval $[0,1]$, belonging respectively to $\mathcal{D}(B(0,\frac{4}{3}))$ and $\mathcal{D}(\mathcal{C})$, and such that $$\forall\xi\in\mathbb{R}^{d},~\chi(\xi)+\sum_{j\geq0}\varphi(2^{-j}\xi)=1,$$
$$\forall\xi\in\mathbb{R}^{d}\backslash\{0\},~\sum_{j\in\mathbb{Z}}\varphi(2^{-j}\xi)=1,$$
$$|j-j'|\geq2\Rightarrow Supp ~\varphi(2^{-j}\xi)\cap Supp ~\varphi(2^{-j'}\xi)=\emptyset,$$
$$j\geq1\Rightarrow Supp ~\chi(\xi)\cap Supp ~\varphi(2^{-j'}\xi)=\emptyset.$$
  Define the set $\widetilde{\mathcal{C}}=B(0,\frac{2}{3})+\mathcal{C}$. And we have
$$|j-j'|\geq5\Rightarrow 2^{j'}\widetilde{\mathcal{C}}\cap 2^{j}\mathcal{C}=\emptyset.$$
Further, we have $$\forall\xi\in\mathbb{R}^{d},~\frac{1}{2}\leq\chi^{2}(\xi)+\sum_{j\geq0}\varphi^{2}(2^{-j}\xi)\leq1,$$
$$\forall\xi\in\mathbb{R}^{d},~\frac{1}{2}\leq\sum_{j\in\mathbb{Z}}\varphi^{2}(2^{-j}\xi)\leq1.$$
   Denote $\mathcal{F}$ by the Fourier transform and $\mathcal{F}^{-1}$ by its inverse. From now on, we write $h=\mathcal{F}^{-1}\varphi$ and $\widetilde{h}=\mathcal{F}^{-1}\chi$.
The nonhomogeneous dyadic blocks $\Delta_{j}$ are defined by
$$\Delta_{j}u=0~~~ if~~~ j\leq-2,~~~\Delta_{-1}u=\chi(D)u=\int_{\mathbb{R}^{d}}\widetilde{h}(y)u(x-y)dy,$$
$$ and,~~~\Delta_{j}u=\varphi(2^{-j}D)u=2^{jd}\int_{\mathbb{R}^{d}}h(2^{j}y)u(x-y)dy ~~~if~~ j\geq0,$$
$$S_{j}u=\sum_{j'\leq j-1}\Delta_{j'}u.$$
    And the homogeneous dyadic blocks $\dot{\Delta}_{j}$ are defined by
$$ \dot{\Delta}_{j}u=\varphi(2^{-j}D)u=2^{jd}\int_{\mathbb{R}^{d}}h(2^{j}y)u(x-y)dy, $$
$$\dot{S}_{j}u=\chi(2^{-j}D)u=\int_{\mathbb{R}^{d}}\widetilde{h}(2^{j}y)u(x-y)dy.$$
   The homogeneous and nonhomogeneous Besov spaces are denote by $\dot{B}^{s}_{p,r}$  and $B^{s}_{p,r}$
$$\dot{B}^{s}_{p,r}=\big\{u\in S'_{h}\big{|}\|u\|_{\dot{B}^{s}_{p,r}}=(\sum_{j\in\mathbb{Z}}2^{rjs}\|\dot{\Delta}_{j}u\|^{r}_{L^{p}})^{\frac{1}{r}}<\infty\big\},$$
$$B^{s}_{p,r}=\big\{u\in S'\big{|}\|u\|_{B^{s}_{p,r}}=(\sum_{j\geq-1}2^{rjs}\|\Delta_{j}u\|^{r}_{L^{p}})^{\frac{1}{r}}<\infty\big\}.$$
   Also we denote $C_{T}(B^{s}_{p,r})$ by $C([0,T];B^{s}_{p,r})$  and $L^{\rho}_{T}(B^{s}_{p,r})$ by $L^{\rho}([0,T];B^{s}_{p,r})$ respectively. Moreover we use the spaces $\widetilde{L}^{\rho}_{T}(B^{s}_{p,r})$ and $B^{s}_{p,r}(\mathcal{L}^{q})$
$$\widetilde{L}^{\rho}_{T}(B^{s}_{p,r})=\big\{u\in S'\big{|}\|u\|_{B^{s}_{p,r}}=(\sum_{j\geq-1}2^{rjs}\|\Delta_{j}u\|^{r}_{L^{\rho}_{t}(L^{p})})^{\frac{1}{r}}<\infty\big\},$$
$$B^{s}_{p,r}(\mathcal{L}^{q})=\big\{\phi\in S'\big{|}\|\phi\|_{B^{s}_{p,r}(\mathcal{L}^{q})}=(\sum_{j\geq-1}2^{rjs}\|\Delta_{j}\phi\|^{r}_{L^{p}_{x}(\mathcal{L}^{q})})^{\frac{1}{r}}<\infty\big\}.$$
   The dyadic blocks $\Delta_{j}$ are related to the variable $x$ and independent with variable $R$. It means that  $$\mathcal{F}(\psi(t,x,R))=\int_{\mathbb{R}^{d}}\psi(t,x,R)e^{-i\xi\cdot x}dx.$$
   Next we define a special space $E^s_{p,r}$ which is useful in this paper,
$$
E^s_{p,r}(T)=\bigg\{\psi:\|\psi\|_{E^{s}_{p,r}(T)}=\bigg(\sum_{j\geq -1}\bigg[2^{js}\bigg(\displaystyle\int^{T}_{0}\displaystyle\int_{\mathbb{R}^{d}\times B}\bigg|\nabla_{R}\bigg(\Delta_{j}\displaystyle\frac{\psi}{\psi_{\infty}}\bigg)^{\frac{p}{2}}\bigg|^{2}\psi_{\infty}dxdRdt\bigg)^{\frac{1}{p}}\bigg]^{r}\bigg)^{\frac{1}{r}}<\infty\bigg\}.$$
   Now we state our main results.
\begin{theo}
Assume that $s>\frac{d}{p}+1,~2\leq p<\infty,~r\geq p$ and $u_{0}\in B^{s}_{p,r}(\mathbb{R}^{d}),~\psi_{0}\in B^{s}_{p,r}(\mathbb{R}^{d};\mathcal{L}^{p})$. Then there exist some $T^{*}>0$ and a unique solution $(u,\psi)$ of $(1.2)$ in
$$C([0,T^{*});B^{s}_{p,r})\times C([0,T^{*});B^{s}_{p,r}(\mathbb{R}^{d};\mathcal{L}^{p})),~~~~~if~r<\infty,$$
$$C_w([0,T^{*});B^{s}_{p,\infty})\times C_w([0,T^{*});B^{s}_{p,\infty}(\mathbb{R}^{d};\mathcal{L}^{p})).$$
Moreover $u\in L^{2}_{T^{*}}(B^{s+1}_{p,r})$ and $\psi \in E^s_{p,r}(T)$.
\end{theo}

\begin{theo}
Let $(u_0,\psi_0)$ be as in Theorem 2.1, and let $T^*$ be the lifespan of the solution to (1.2). If $T^*<\infty$, then we have that
$$\int^{T^*}_0\|u\|^2_{L^\infty}=\infty,$$
$$\|\psi\|^2_{L^\infty_{T^{*}}(B^{s}_{p,r}(\mathcal{L}^{p}))}+\|\psi\|^2_{E^s_{p,r}(T*)}=\infty.$$
\end{theo}

\begin{theo}
Under the assumption of Theorem 2.1  and assume that $\int_{B}\psi_{0}dR=1 ~~a.e. ~~in~~x$. If there exists a constant $c_{0}$ such that $$\|u_{0}\|^{2}_{B^{s}_{p,r}}+\|\psi_{0}-\psi_{\infty}\|^{2}_{B^{s}_{p,r}(\mathcal{L}^{p})}\leq c_{0},$$
then the solution constructed in Theorem 2.1 is global. Moreover, there exist a constant $M$ such that
$$\|u(t)\|^{2}_{B^{s}_{p,r}}+\|\psi(t)-\psi_{\infty}\|^{2}_{B^{s}_{p,r}(\mathcal{L}^{p})}\leq Mc_{0}.$$
\end{theo}

\begin{rema}
Thanks to $B^s_{2,2}=H^s$ for $s>0$, if we take $p=2, r=2$ in Theorem 2.1 and Theorem 2.3, then our theorems cover the recent results obtained by Masmoudi in \cite{Masmoudi.W}.
\end{rema}

\section{Preliminaries}
In this section we introduce some useful lemmas which will be used in the sequel. For more details, one can refer to Section 2 in \cite{B.C.D}.
\subsection{Propositions of Besov spaces}
Firstly we introduce the Bernstein inequalities.
\begin{lemm}
\cite{B.C.D} Let $\mathcal{C}$ be an annulus and $B$ a ball. A constant $C$ exists such that
for any nonnegative integer $k$, any couple $(p, q)$ in $[1,\infty]^{2}$ with $q \geq p\geq 1$, and
any function $u$ of $L^{p}$, we have
$$Supp ~\widehat{u} \subseteq \lambda B \Rightarrow \|D^{k}u\|_{L^{q}}\triangleq \sup_{|\alpha|\leq k} \|\partial ^{\alpha}u\|_{L^{q}}\leq C^{k+1}\lambda^{k+d(\frac{1}{p}-\frac{1}{q})}\|u\|_{L^{p}},$$
$$Supp ~\widehat{u} \subseteq \lambda\mathcal{C} \Rightarrow C^{-k-1}\lambda^{k}\|u\|_{L^{p}} \leq\|D^{k}u\|_{L^{p}}\leq C^{k+1}\lambda^{k}\|u\|_{L^{p}}.$$
\end{lemm}

\begin{lemm}
\cite{B.C.D} Let $\mathcal{C}$  be an annulus. Positive constants $c$ and $C$ exist such that
for any $p$ in $[1,\infty]$ and any couple $(t,\lambda )$ of positive real numbers, we have
$$Supp ~\widehat{u} \subseteq \lambda\mathcal{C} \Rightarrow  \|e^{t\triangle}u\|_{L^{p}}\leq Ce^{-ct\lambda^{2}}\|u\|_{L^{p}}.$$
\end{lemm}

The following proposition is about the embedding for Besov spaces.
\begin{prop}
\cite{B.C.D} Let $1\leq p_{1} \leq p_{2} \leq \infty$ and $1\leq r_{1} \leq r_{2} \leq \infty$, and let $s$ be a real number. Then we have
$$B^{s}_{p_{1},r_{1}}\hookrightarrow B^{s-d(\frac{1}{p_{1}}-\frac{1}{p_{2}})}_{p_{2},r_{2}}.$$
If $s>\frac{d}{p}~or ~s=\frac{d}{p},~r=1$, we then have $$B^{s}_{p,r}\hookrightarrow L^{\infty}.$$
\end{prop}

\begin{prop}
\cite{B.C.D} The set $B^{s}_{p,r}$ is a Banach space and satisfies the Fatou property, namely, if $(u_{n})_{n\in\mathbb{N}}$ is a bounded sequence of $B^{s}_{p,r}$, then an element $u$ of $B^{s}_{p,r}$ and a subsequence $(u_{\psi(n)})$ exist such that
$$\lim_{n\rightarrow\infty}u_{\psi(n)}=u ~~~in~~~ S' ~~ and ~~\|u\|_{B^{s}_{p,r}}\leq\ C\liminf_{n\rightarrow\infty}\|u_{\psi(n)}\|_{B^{s}_{p,r}}.$$

\end{prop}

Next we introduce the Bony decomposition.
\begin{defi}
\cite{B.C.D} The nonhomogeneous paraproduct of $v$ and $u$ is defined by
$$T_{u}v\triangleq\sum_{j}S_{j-1}u\Delta_{j}v.$$
The nonhomogeneous remainder of $v$ and $u$ is defined by
$$R(u,v)\triangleq\sum_{|k-j|\leq1}\Delta_{k}u\Delta_{j}v.$$
We have the following  Bony decomposition
$$uv=T_{u}v+R(u,v)+T_{v}u.$$
\end{defi}

\begin{prop}
\cite{B.C.D} For any couple of real numbers $(s, t)$ with $t$ negative and any $(p, r_{1}, r_{2})$
in $[1,\infty]^{3}$, there exists a constant $C$ such that:
$$\|T_{u}v\|_{B^{s}_{p,r}}\leq C\|u\|_{L^{\infty}}\|D^{k}v\|_{B^{s-k}_{p,r}},$$
$$\|T_{u}v\|_{B^{s+t}_{p,r}}\leq C\|u\|_{B^{t}_{p,r_{1}}}\|D^{k}v\|_{B^{s-k}_{p,r_{2}}},$$
where $r=\min\{1,\frac{1}{r_{1}}+\frac{1}{r_{2}}\}.$
\end{prop}

\begin{prop}
\cite{B.C.D} A constant $C$ exists which satisfies the following inequalities.
Let $(s_{1}, s_{2})$ be in $\mathbb{R}^2$ and $(p_{1}, p_{2}, r_{1}, r_{2})$ be in $[1,\infty]^{4}$. Assume that
$$\frac{1}{p}=\frac{1}{p_{1}}+\frac{1}{p_{2}}~~and~~\frac{1}{r}=\frac{1}{r_{1}}+\frac{1}{r_{2}} .$$
If $s_{1}+s_{2}>0$, then we have, for any $(u,v)$ in $B^{s_{1}}_{p_{1},r_{1}}\times B^{s_{2}}_{p_{2},r_{2}}$,
$$\|R(u,v)\|_{B^{s_{1}+s_{2}}_{p,r}}\leq \frac{C^{s_{1}+s_{2}+1}}{s_{1}+s_{2}}\|u\|_{B^{s_{1}}_{p_{1},r_{1}}}\|v\|_{B^{s_{2}}_{p_{2},r_{2}}}.$$
If $r=1$ and $s_{1}+s_{2}=0$, then we have, for any $(u,v)$ in $B^{s_{1}}_{p_{1},r_{1}}\times B^{s_{2}}_{p_{2},r_{2}}$,
$$\|R(u,v)\|_{B^{0}_{p,\infty}}\leq C\|u\|_{B^{s_{1}}_{p_{1},r_{1}}}\|v\|_{B^{s_{2}}_{p_{2},r_{2}}}.$$
\end{prop}

\begin{coro}
\cite{B.C.D} For any positive real number $s$ and any $(p, r)$ in $[1,\infty]^{2}$, the
space $L^{\infty}\cap B^{s}_{p,r}$ is an algebra, and a constant $C$ exists such that
$$\|uv\|_{B^{s}_{p,r}}\leq C(\|u\|_{L^{\infty}}\|v\|_{B^{s}_{p,r}}+\|u\|_{B^{s}_{p,r}}\|v\|_{L^{\infty}}).$$
If $s>\frac{d}{p}$ or $s=\frac{d}{p},~r=1$, we have
$$\|uv\|_{B^{s}_{p,r}}\leq C\|u\|_{B^{s}_{p,r}}\|v\|_{B^{s}_{p,r}}.$$
\end{coro}

\subsection{Propositions of $\mathcal{L}^{p}$  and $B^{s}_{p,r}(\mathcal{L}^{p})$ spaces}
Now we introduce some propositions about $\mathcal{L}^{p}$ and $B^{s}_{p,r}(\mathcal{L}^{p})$.

\begin{prop}
For any $1\leq p \leq \infty$. We have $\mathcal{L}^{p}(B)\hookrightarrow L^{p}(B)$.
\begin{proof}
By the definition we have  $$\psi_{\infty}=\frac{e^{-\mathcal{U}(R)}}{\int_{B}e^{-\mathcal{U}(R)}dR}=\frac{(1-|R|^{2})^{k}}{\int_{B}(1-|R|^{2})^{k}dR}.$$ Since $|R|\leq 1$, it follows that $\psi_{\infty}^{p}\leq C \psi_{\infty}$, where $C$ is a constant independent of $R$. Then
$$\int_{B}|\psi|^{p}=\int_{B}\psi_{\infty}^{p}\frac{|\psi|^{p}}{\psi_{\infty}^{p}}\leq C \int_{B}\psi_{\infty}\frac{|\psi|^{p}}{\psi_{\infty}^{p}}. $$
\end{proof}
\end{prop}

\begin{rema}
If $1\leq p<\infty$, then $C^{\infty}_{0}(B)$ is dense in $\mathcal{L}^{p}(B)$.
\end{rema}

\begin{prop}
If $1<p<\infty$ then $\mathcal{L}^{p}$ is a reflexive Banach space.
\begin{proof}
By the definition of $\mathcal{L}^{p}$, we have
$$\mathcal{L}^{p}=L^{p}(\frac{dR}{\psi^{p-1}_{\infty}}).$$
Since $\frac{dR}{\psi^{p-1}_{\infty}}$ is a $\sigma$ finitely additive measure,  it follows that
$$(\mathcal{L}^{p})^{*}=(L^{p}(\frac{dR}{\psi^{p-1}_{\infty}}))^{*}=L^{p'}(\frac{dR}{\psi^{p-1}_{\infty}}),$$
and $$(L^{p'}(\frac{dR}{\psi^{p-1}_{\infty}}))^{*}=L^{p}(\frac{dR}{\psi^{p-1}_{\infty}}).$$
So $\mathcal{L}^{p}$ is a reflexive Banach space.
\end{proof}
\end{prop}

Next we will introduce an inequality which is used to estimate the stress tensor $div~\tau$, and we set $x=1-|R|$.
\begin{lemm}
\cite{Masmoudi.W} For all $\varepsilon>0$, there exists a constant $C_{\varepsilon}$ such that
$$(\int_{B}\frac{|\psi|}{x}dR)^{2}\leq\varepsilon\int_{B}\psi_{\infty}|\nabla_{R}\frac{\psi}{\psi_{\infty}}|^{2}dR+C_{\varepsilon}\int_{B}\frac{|\psi|^{2}}{\psi_{\infty}}dR.$$
\end{lemm}

\begin{coro}
\cite{Kreml} Assume that $2\leq p<\infty$, for all $\varepsilon>0$, there exists a constant $C_{\varepsilon}$ such that
$$(\int_{B}\frac{|\psi|}{x}dR)^{p}\leq\varepsilon\int_{B}\psi_{\infty}|\nabla_{R}(\frac{\psi}{\psi_{\infty}})^{\frac{p}{2}}|^{2}dR+C_{\varepsilon}\int_{B}|\frac{\psi}{\psi_{\infty}}|^{p}\psi_{\infty}dR.$$
\begin{proof}
By a direct calculation and the H\"{o}lder inequality, we have
\begin{align}
\nonumber (\int_{B}\frac{|\psi|}{x}dR)^{p}&=(\int_{B}x^{\frac{2}{p}(k-1)-k}|\psi|x^{\frac{(p-2)(k-1)}{(p)}} dR)^{p}\\
\nonumber &\leq(\int_{B}x^{(k-1)-\frac{kp}{2}}|\psi|^{\frac{p}{2}}dR)^{2}(\int_{B}x^{k-1}dR)^{p-2}\\
\nonumber &\leq C(\int_{B}x^{(k-1)-\frac{kp}{2}}|\psi|^{\frac{p}{2}}dR)^{2} \\
&=C (\int_{B}\frac{x^{k-\frac{kp}{2}}|\psi|^{\frac{p}{2}}}{x}dR)^{2}.
\end{align}
So by Lemma (3.12) we obtain
\begin{align*}
 (3.1)&\leq\varepsilon\int_{B}\psi_{\infty}|\nabla_{R}\frac{(\psi_{\infty})^{1-\frac{p}{2}}(\psi)^{\frac{p}{2}}}{\psi_{\infty}}|^{2}dR+C_{\varepsilon}\int_{B}|\frac{(\psi_{\infty})^{1-\frac{p}{2}}(\psi)^{\frac{p}{2}}}{\psi_{\infty}}|^{2}dR\\
 &=\varepsilon\int_{B}\psi_{\infty}|\nabla_{R}(\frac{\psi}{\psi_{\infty}})^{\frac{p}{2}}|^{2}dR+C_{\varepsilon}\int_{B}|\frac{\psi}{\psi_{\infty}}|^{p}\psi_{\infty}.
\end{align*}
\end{proof}
\end{coro}

Next we will introduce some propositions about the space $B^{s}_{p,r}(\mathcal{L}^{q})$.
The following proposition is similar to Proposition (3.4). For more details, one can refer to Section 2.3 in \cite{B.C.D}.
\begin{prop}
The set $B^{s}_{p,r}(\mathcal{L}^{q})$ is a Banach space. Moreover if $1< p,~q<\infty $, the set $B^{s}_{p,r}(\mathcal{L}^{q})$ satisfies the Fatou property, namely, if $(\psi_{n})_{n\in\mathbb{N}}$ is a bounded sequence of $B^{s}_{p,r}(\mathcal{L}^{q})$, then an element $\psi$ of $B^{s}_{p,r}(\mathcal{L}^{q})$ and a subsequence $(\psi_{n_{k}})$ exist such that
$$\lim_{k\rightarrow\infty}\psi_{n_{k}}=\psi ~~~in~~~ S' ~~ and ~~\|\psi\|_{B^{s}_{p,r}(\mathcal{L}^{q})}\leq\ C\liminf_{k\rightarrow\infty}\|\psi_{n_{k}}\|_{B^{s}_{p,r}(\mathcal{L}^{q})} .$$
\begin{proof}
Assume that the sequence $(\psi_{n})$ is bounded in $B^{s}_{p,r}(\mathcal{L}^{q})$. Then for any $j\geq-1$, the sequence $(\Delta_{j}\psi_{n})$ is bounded in $L^{p}_{x}(\mathcal{L}^{q})$. Because $\mathcal{L}^{q}$ is a reflexive Banach space. Cantor’s diagonal process thus supplies a subsequence $(\Delta_{j}\psi_{n_{k}})$ and a sequence $(\widetilde{\psi}_{j})$ of $C^{\infty}(\mathbb{R}^{d};\mathcal{L}^{q})$ functions with Fourier transform supported in $2^{j}\mathcal{C}$ such that, for any $j\in\mathbb{Z}$, $\phi\in S$,
\begin{align*}
\lim_{k\rightarrow\infty}\langle\Delta_{j}\psi_{n_{k}},~\phi\rangle=\langle\widetilde{\psi}_{j},~\phi\rangle~~and~~\|\widetilde{\psi_{j}}\|_{L^{p}(\mathcal{L}^{q})}\leq \liminf_{n\rightarrow\infty}\|\psi_{n}\|_{L^{p}(\mathcal{L}^{q})}.
\end{align*}
Now, the sequence $((2^{js}\|\Delta_{j}\psi_{n_{k}}\|_{L^{p}(\mathcal{L}^{q})})_{j})_{k}$ is bounded in $l^{r}$. Hence, there exists an element $\widetilde{c}_{j}$ of $l^{r}$ such that (up to an extraction)
 for any sequence $(d_{j})$ of nonnegative real numbers different from $0$
for only a finite number of indices $j$,
$$\lim_{k\rightarrow\infty}\sum_{j\geq-1}2^{js}\|\Delta_{j}\psi_{n_{k}}\|_{L^{p}_{x}(\mathcal{L}^{q})}d_{j}=\sum_{j\geq-1}\widetilde{c}_{j}d_{j}~~and ~~\|\widetilde{c}_{j}\|_{l^{r}}\leq\liminf_{k\rightarrow\infty}\|\psi_{n_{k}}\|_{B^{s}_{p,r}(\mathcal{L}^{q})}.$$
Passing to the limit in the sum gives that $(2^{js}\|\widetilde{\psi}_{j}\|_{L^{p}_{x}(\mathcal{L}^{q})})_{j}$ belongs to $l^{r}$. The Fourier transform of $(\widetilde{\psi}_{j})$ is supported in $2^{j}\mathcal{C}$. So the series $\sum_{j\geq-1}\widetilde{\psi}_{j}$ converges to some $\psi$ in $S'$. For all the $N$ and $\phi\in S$ we have
$$\big{<}\sum^{N}_{j=-1}\Delta_{j}\psi,\phi\big{>}=\big{<}\sum^{N}_{j=-1}\sum_{|j'-j|\leq1}\Delta_{j}\widetilde{\psi}_{j},\phi\big{>},$$
then we obtain
$$\sum^{N}_{j=-1}\Delta_{j}\psi=\lim_{k\rightarrow\infty}\sum^{N}_{j=-1}\Delta_{j}\psi_{n_{k}} ~~in~~S' ,$$
and $(Id-S_{N})\psi_{n_{k}}$ tends to $0$ in $B^{s}_{p,r}$. So $\psi$ is indeed the limit of $\psi_{n_{k}}$ in $S'$, which completes the proof of the Fatou property.

We will now check that $B^{s}_{p,r}(\mathcal{L}^{q})$ is complete. Consider a Cauchy sequence $\psi_{n}$. This sequence is of course bounded, so there exist some $\psi$ and a subsequence $\psi_{n_{k}}$, such that $\psi_{n_{k}}$ converges to $\psi$ in $S'$. For any $\varepsilon>0$ there exist a $N_{\varepsilon}$ such that
   $$n_{k}\geq m_{k}\geq N_{\varepsilon}\Rightarrow \|\psi_{n_{k}}-\psi_{m_{k}}\|_{B^{s}_{p,r}(\mathcal{L}^{q})}<\varepsilon.$$
   The Fatou property ensures that
   $$n_{k}\geq N_{\varepsilon},~\|\psi_{n_{k}}-\psi\|_{B^{s}_{p,r}(\mathcal{L}^{q})}<\varepsilon.$$
   Hence $\psi_{n_{k}}$ tends to $\psi$ in $B^{s}_{p,r}(\mathcal{L}^{q})$, as $k$ goes to infinity.

\end{proof}
\end{prop}

The next lemma is very useful to deal with the product of $u(t,x)$ and $\psi(t,x,R)$.
\begin{lemm}
For any positive real number $s$ and any $(p, r)$ in $[1,\infty]^{2}$, if $u=u(t,x)\in L^{\infty}\cap B^{s}_{p,r}$ and $\psi=\psi(t,x,R)\in L^{\infty}_{x}(\mathcal{L}^{p})\cap B^{s}_{p,r}(\mathcal{L}^{p})$, then a constant $C$ exists such that
$$\|u\psi\|_{B^{s}_{p,r}(\mathcal{L}^{p})}\leq C(\|u\|_{L^{\infty}}\|\psi\|_{B^{s}_{p,r}(\mathcal{L}^{p})}+\|u\|_{B^{s}_{p,r}}\|\psi\|_{L^{\infty}_{x}(\mathcal{L}^{p})}).$$
If $s>\frac{d}{p}$ or $s=\frac{d}{p},~r=1$, we have
$$\|u\psi\|_{B^{s}_{p,r}(\mathcal{L}^{p})}\leq C\|u\|_{B^{s}_{p,r}}\|\psi\|_{B^{s}_{p,r}(\mathcal{L}^{p})}.$$
\begin{proof}
 We can write $$u\psi=T_{u}\psi+T_{\psi}u+R(u,\psi).$$
 Firstly, we consider the term $T_{u}\psi$, by definition and the proposition about $\Delta_{j}~and~S_{j}$ we deduce that
$$\Delta_{j}T_{u}\psi=\Delta_{j}\sum_{|j-j'|\leq2}S_{j'-1}u\Delta_{j'}\psi.$$
So using H\"{o}lder's inequality, we infer that
\begin{align*}
2^{js}\|\Delta_{j}T_{u}\psi\|_{L^{p}_{x}(\mathcal{L}^{p})}&\leq2^{js}\sum_{|j-j'|\leq2}\|S_{j'-1}u\Delta_{j'}\psi\|_{L^{p}_{x}(\mathcal{L}^{p})}\\
&=\sum_{|j-j'|\leq2}2^{js}\bigg(\int_{\mathbb{R}^{d}\times B}\bigg|S_{j'-1}u(t,x)\frac{\Delta_{j'}\psi(t,x,R)}{\psi_{\infty}}\bigg|^{p}\psi_{\infty}dxdR\bigg)^{\frac{1}{p}}\\
&\leq C\sum_{|j-j'|\leq2}2^{js}\|S_{j'-1}u(t,x)\|_{L^{\infty}}\bigg(\int_{\mathbb{R}^{d}\times B}\bigg|\frac{\Delta_{j'}\psi(t,x,R)}{\psi_{\infty}}\bigg|^{p}\psi_{\infty}dxdR\bigg)^{\frac{1}{p}}\\
&\leq C\sum_{|j-j'|\leq2}2^{(j-j')s}\|u(t,x)\|_{L^{\infty}}2^{j's}\|\Delta_{j'}\psi\|_{L^{p}_{x}(\mathcal{L}^{p})}\\
&\leq C c_{j}\|u\|_{L^{\infty}}\|\psi\|_{B^{s}_{p,r}(\mathcal{L}^{p})}.
\end{align*}
Where $c_{j}$ denotes an element of the unit sphere of $l^{r}$. Taking $l^{r}$-norm for both sides of the above inequality, we obtain
\begin{align}
\|T_{u}\psi\|_{B^{s}_{p,r}(\mathcal{L}^{p})}\leq C \|u\|_{L^{\infty}}\|\psi\|_{B^{s}_{p,r}(\mathcal{L}^{p})}.
\end{align}
Next, we consider the second term $T_{\psi}u$. Similarly we deduce that
$$\Delta_{j}T_{\psi}u=\Delta_{j}\sum_{|j-j'|\leq2}S_{j'-1}\psi\Delta_{j'}u.$$
So using Fubini's theorem and H\"{o}lder's inequality, we infer that
\begin{align*}
2^{js}\|\Delta_{j}T_{\psi}u\|_{L^{p}_{x}(\mathcal{L}^{p})}&\leq2^{js}\sum_{|j-j'|\leq2}\|S_{j'-1}\psi\Delta_{j'}u\|_{L^{p}_{x}(\mathcal{L}^{p})}\\
&=\sum_{|j-j'|\leq2}2^{js}\bigg(\int_{\mathbb{R}^{d}\times B}\bigg|\Delta_{j'}u(t,x)\frac{S_{j'-1}\psi(t,x,R)}{\psi_{\infty}}\bigg|^{p}\psi_{\infty}dxdR\bigg)^{\frac{1}{p}}\\
&=\sum_{|j-j'|\leq2}2^{js}(\int_{\mathbb{R}^{d} }|\Delta_{j'}u(t,x)|^{p}\int_{B}\bigg|\frac{S_{j'-1}\psi(t,x,R)}{\psi_{\infty}}\bigg|^{p}\psi_{\infty}dRdx)^{\frac{1}{p}}\\
&\leq C\sum_{|j-j'|\leq2}2^{js}\|\Delta_{j'}u(t,x)\|_{L^{p}}\sup_{x}\bigg(\int_{ B}\bigg|\frac{S_{j'-1}\psi(t,x,R)}{\psi_{\infty}}\bigg|^{p}\psi_{\infty}dxdR\bigg)^{\frac{1}{p}}\\
&\leq C\sum_{|j-j'|\leq2}2^{(j-j')s}2^{j's}\|\Delta_{j'}u(t,x)\|_{L^{p}}\|\psi\|_{L^{\infty}_{x}(\mathcal{L}^{p})}\\
&\leq C c_{j}\|u\|_{B^{s}_{p,r}}\|\psi\|_{L^{\infty}(\mathcal{L}^{p})}.
\end{align*}
Where $c_{j}$ denotes an element of the unit sphere of $l^{r}$. Taking $l^{r}$-norm for both sides of the above inequality, we get
\begin{align}
\|T_{\psi}u\|_{B^{s}_{p,r}(\mathcal{L}^{p})}\leq C \|u\|_{B^{s}_{p,r}}\|\psi\|_{L^{\infty}(\mathcal{L}^{p})}.
\end{align}
Finally, we consider the last term $R(u,\psi)$. By definition, we can write
$$R(u,\psi)=\sum_{j'}R_{j'},~~with~~~ R_{j'}=\sum_{|k|\leq1}\Delta_{j'-k}u(t,x)\Delta_{j'}\psi(t,x,R).$$
By the construction of the dyadic partition of unity, there
exists an integer $N_{0}$ such that
$$j>j'+N_{0}\Rightarrow  \Delta_{j}R_{j'}=0.$$
From this we deduce that
$$\Delta_{j}R(u,\psi)=\Delta_{j}\sum_{j'\geq j-N_{0}}\sum_{|k|\leq1}\Delta_{j'-k}u\Delta_{j'}\psi.$$
So using H\"{o}lder's inequality and due to $s>0$, we infer that
\begin{align*}
2^{js}\|\Delta_{j}R(u,\psi)\|_{L^{p}_{x}(\mathcal{L}^{p})}&\leq2^{js}\sum_{j'\geq j-N_{0}}\sum_{|k|\leq1}\|\Delta_{j'-k}u\Delta_{j'}\psi\|_{L^{p}_{x}(\mathcal{L}^{p})}\\
&=\sum_{j'\geq j-N_{0}}\sum_{|k|\leq1}2^{js}\bigg(\int_{\mathbb{R}^{d}\times B}\bigg|\Delta_{j'-k}u(t,x)\frac{\Delta_{j'}\psi(t,x,R)}{\psi_{\infty}}\bigg|^{p}\psi_{\infty}dxdR\bigg)^{\frac{1}{p}}\\
&\leq C\sum_{j'\geq j-N_{0}}2^{js}\|u(t,x)\|_{L^{\infty}}\bigg(\int_{\mathbb{R}^{d}\times B}\bigg|\frac{\Delta_{j'}\psi(t,x,R)}{\psi_{\infty}}\bigg|^{p}\psi_{\infty}dxdR\bigg)^{\frac{1}{p}}\\
&\leq C\sum_{j'\geq j-N_{0}}2^{(j-j')s}\|u(t,x)\|_{L^{\infty}}2^{j's}\|\Delta_{j'}\psi\|_{L^{p}_{x}(\mathcal{L}^{p})}\\
&\leq C c_{j}\|u\|_{L^{\infty}}\|\psi\|_{B^{s}_{p,r}(\mathcal{L}^{p})}.
\end{align*}
Where $c_{j}$ denotes an element of the unit sphere of $l^{r}$. Taking $l^{r}$-norm for both sides of the above inequality, we obtain
\begin{align}
\|R(u,\psi)\|_{B^{s}_{p,r}(\mathcal{L}^{p})}\leq C \|u\|_{L^{\infty}}\|\psi\|_{B^{s}_{p,r}(\mathcal{L}^{p})}.
\end{align}
Combining (3.2), (3.3) and (3.4), we complete the proof.
\end{proof}
\end{lemm}

\subsection{Commutator estimates}
This section is devoted to various commutator estimates which enable us to establish a priori estimates. The proof is similar to the commutator estimates for Besov spaces. For more details, one can refer to Section 2.10 in \cite{B.C.D}.
\begin{lemm}
Let $\theta$ be a $C^{1}$ function on $\mathbb{R}^{d}$ such that $|\cdot|\mathcal{F}^{-1}(\theta)\in L^{1}$. There exists a constant $C$ such that for any Lipschitz function $u(x)$ and any function $\psi(x,R)$ in $L^{p}_{x}(\mathcal{L}^{p})$, we have, for any positive $\lambda$,
$$\|[\theta(\lambda^{-1}D),u]\psi\|_{L^{p}_{x}(\mathcal{L}^{p})}\leq C\lambda^{-1}\|\nabla u\|_{L^{\infty}}\|\psi\|_{L^{p}_{x}(\mathcal{L}^{p})}.$$
\begin{proof}
Indeed,
$$([\theta(\lambda^{-1}D),u]\psi)=\lambda^{d}\int_{\mathbb{R}^{d}}k(\lambda(x-y))(u(y)-u(x))\psi(y,R)dy~~~with~~k=\mathcal{F}^{-1}\theta.$$
Let $k_{1}(z)=|z||k(z)|$. From the first order Taylor formula, we deduce that
$$([\theta(\lambda^{-1}D),u]\psi)(x,R)\leq\lambda^{-1}\int_{[0,1]\times\mathbb{R}^{d}}\lambda^{d}k_{1}(\lambda z)(\nabla |u(x-\tau z))||\psi(x-z,R)|dz.$$
Taking the $L^{p}_{x}(\mathcal{L}^{p})$-norm of the above inequality and using  H\"{o}lder's inequality, we infer that
\begin{align*}
\|[\theta(\lambda^{-1}D),u]\psi\|_{L^{p}_{x}(\mathcal{L}^{p})}&\leq\lambda^{-1}\int_{[0,1]\times\mathbb{R}^{d}}\lambda^{d}k_{1}(\lambda z)(\nabla \|u(x-\tau z))\psi(x-z,R)\|_{L^{p}_{x}(\mathcal{L}^{p})}dz\\
&\leq C\lambda^{-1}\|k_{1}\|_{L^{1}}\|\nabla u(x)\psi(x,R)\|_{L^{p}_{x}(\mathcal{L}^{p})} \\
&\leq C\lambda^{-1}\int_{\mathbb{R}^{d}\times B}\bigg|\nabla u(x)\frac{\psi(x,R)}{\psi_{\infty}}\bigg|^{p}\psi_{\infty}dxdR)^{\frac{1}{p}}\\
&\leq C\lambda^{-1}\|\nabla u\|_{L^{\infty}}\int_{\mathbb{R}^{d}\times B}\bigg|\frac{\psi(x,R)}{\psi_{\infty}}\bigg|^{p}\psi_{\infty}dxdR)^{\frac{1}{p}}\\
&\leq C\lambda^{-1}\|\nabla u\|_{L^{\infty}}\|\psi\|_{L^{p}_{x}(\mathcal{L}^{p})}.
\end{align*}
\end{proof}
\end{lemm}

\begin{lemm}
Let $s>0$ and $(p,r)\in[1,\infty]^{2}$, and let $u$ be a vector field on $\mathbb{R}^{d}$ with $\nabla u\in L^{\infty}\cap B^{s-1}_{p,r}.$
\\Define $R_{j}=[u\cdot\nabla,\Delta_{j}]\psi(t,x,R)$. There exists a constant $C$, depending on $p,s,r~and~d$, such that
$$\big\|\big(2^{js}\|R_{j}\|_{L^{p}_{x}(\mathcal{L}^{p})}\big)_{j}\big\|_{l^{r}}\leq C (\|\nabla u\|_{L^{\infty}}\|\psi\|_{B^{s}_{p,r}(\mathcal{L}^{p})}+\|\nabla_{x}\psi\|_{L^{\infty}_{x}(\mathcal{L}^{p})}\|\nabla u\|_{B^{s-1}_{p,r}}).$$
If $s>1+\frac{d}{p}$ or $s=1+\frac{d}{p},~r=1$, we have
$$\big\|\big(2^{js}\|R_{j}\|_{L^{p}_{x}(\mathcal{L}^{p})}\big)_{j}\big\|_{l^{r}}\leq C\|\nabla u\|_{B^{s-1}_{p,r}}\|\psi\|_{B^{s}_{p,r}(\mathcal{L}^{p})}.$$
If $s\leq1+\frac{d}{p}$, we have
$$\big\|\big(2^{js}\|R_{j}\|_{L^{p}_{x}(\mathcal{L}^{p})}\big)_{j}\big\|_{l^{r}}\leq C\|\nabla u\|_{B^{\frac{d}{p}}_{p,\infty}\cap L^{\infty}}\|\psi\|_{B^{s}_{p,r}(\mathcal{L}^{p})}.$$
\begin{proof}
We shall split $u$ into low and high frequencies: $u=S_{0}u+\widetilde{u}$. Obviously, we have
\begin{align}
\|S_{0}\nabla u\|_{L^{p}}\leq C\|\nabla u\|_{L^{p}} ~~and ~~\|\nabla \widetilde{u}\|_{L^{p}}\leq C\|\nabla u\|_{L^{p}}.
\end{align}
Further, as $\widetilde{u}$ is spectrally supported away from the origin, Lemma 3.1 ensures that
\begin{align}
\forall j\geq-1 ,~~ \|\Delta_{j}\nabla \widetilde{u}\|_{L^{p}}\approx2^{j}\|\Delta_{j} \widetilde{u}\|_{L^{p}}.
\end{align}
Using Bony's decomposition, we end up with $R_{j}=\sum^{8}_{i=1}R^{i}_{j}$, where
\begin{align*}
&R^{1}_{j}=[T_{\widetilde{u}^{k}},~~\Delta_{j}]\partial_{k}\psi,~~~~~R^{2}_{j}=T_{\partial_{k}\Delta_{j}\psi}\widetilde{u}^{k},~~~~~R^{3}_{j}=-\Delta_{j}T_{\partial_{k}\psi}\widetilde{u}^{k},~~~~~R^{4}_{j}=\partial_{k}R(\Delta_{j}\psi,\widetilde{u}^{k}),\\
&R^{5}_{j}=-R(\Delta_{j}\psi,~~div~\widetilde{u}),~~~R^{6}_{j}=-\partial_{k}\Delta_{j}R(\psi,\widetilde{u}^{k}),~~~R^{7}_{j}=\Delta_{j}R(\psi,-div~\widetilde{u})),~~~R^{8}_{j}=[S_{0}u,\Delta_{j}]\partial_{k}\psi.
\end{align*}
In the following computations, we denote by $(c_{j})_{j\geq -1}$ a sequence such that $\|c_{j}\|_{l^{r}}\leq 1$.\\
\textit{Bounds for $2^{js}\|R^{1}_{j}\|_{L^{p}_{x}(\mathcal{L}^{p})}$}. By the construction of the dyadic partition of unity, we have
$$R^{1}_{j}=\sum_{|j-j'|\leq 4}[S_{j'-1}\widetilde{u}^{k},\Delta_{j}]\partial_{k}\Delta_{j'}\psi.$$
Hence, according to Lemma 3.16 and (3.5), we have
\begin{align*}
2^{js}\|R^{1}_{j}\|_{L^{p}_{x}(\mathcal{L}^{p})}&\leq C \|\nabla u \|_{L^{\infty}}\sum_{|j-j'|\leq 4}2^{j's}\|\Delta_{j'}\psi\|_{L^{p}_{x}(\mathcal{L}^{p})}\leq C c_{j}\|\nabla u \|_{L^{\infty}}\|\psi\|_{B^{s}_{p,r}(\mathcal{L}^{p})}.
\end{align*}
\textit{Bounds for $2^{js}\|R^{2}_{j}\|_{L^{p}_{x}(\mathcal{L}^{p})}$.} By the construction of the dyadic partition of unity, we obtain
$$R^{2}_{j}=\sum_{j'\geq j-3}S_{j'-1}\partial_{k}\Delta_{j}\psi\Delta_{j'}\widetilde{u}^{k}.$$
Hence, using (3.5), (3.6) and H\"{o}lder's inequality, we deduce that
$$2^{js}\|R^{2}_{j}\|_{L^{p}_{x}(\mathcal{L}^{p})}\leq C c_{j}\|\nabla u \|_{L^{\infty}}\|\psi\|_{B^{s}_{p,r}(\mathcal{L}^{p})}.$$
\textit{Bounds for $2^{js}\|R^{3}_{j}\|_{L^{p}_{x}(\mathcal{L}^{p})}$}. By the definition, we have
$$R^{3}_{j}=\sum_{|j-j'|\leq4}\Delta_{j}(S_{j'-1}\partial_{k}\psi\Delta_{j'}\widetilde{u}^{k}).$$
Hence, using (3.5), (3.6) and H\"{o}lder's inequality, we obtain
$$2^{js}\|R^{3}_{j}\|_{L^{p}_{x}(\mathcal{L}^{p})}\leq C\sum_{|j-j'|\leq4}2^{(j-j')s}\|\nabla S_{j'-1}\psi\|_{L^{\infty}_{x}(\mathcal{L}^{p})}2^{j'(s-1)}\|\Delta_{j'}\nabla u \|_{L^{p}},$$
from which it follows that
$$2^{js}\|R^{3}_{j}\|_{L^{p}_{x}(\mathcal{L}^{p})}\leq C c_{j}\|\nabla u \|_{B^{s-1}_{p,r}}\|\nabla \psi\|_{L^{\infty}_{x}(\mathcal{L}^{p})}.$$
If $s<1+\frac{d}{p}$, we may write $R^{3}_{j}$ as follow:
\begin{align*}
R^{3}_{j}=\sum_{\begin{subarray}-{}
 |j-j'|\leq4 \\
j''<j'-2
\end{subarray}}
\Delta_{j}(\Delta_{j''}\partial_{k}\psi\Delta_{j'}\widetilde{u}^{k}).
\end{align*}
Hence, using (3.5), (3.6) and H\"{o}lder's inequality, we get
\begin{align*}
2^{js}\|R^{3}_{j}\|_{L^{p}_{x}(\mathcal{L}^{p})}&\leq C\sum_{\begin{subarray}-{}
 |j-j'|\leq4 \\
j''<j'-2
\end{subarray}}2^{js}\|\Delta_{j''}\psi\|_{L^{\infty}_{x}(\mathcal{L}^{p})}2^{-j'}\|\Delta_{j'}\nabla u \|_{L^{p}} \\
&\leq C\sum_{\begin{subarray}-{}
 |j-j'|\leq4 \\
j''<j'-2
\end{subarray}}2^{(j-j'')(s-1-\frac{d}{p})}2^{j''s}\|\Delta_{j''}\psi\|_{L^{\infty}_{x}(\mathcal{L}^{p})}2^{j'\frac{d}{p}}\|\Delta_{j'}\nabla u \|_{L^{p}}.
\end{align*}
Since $s<1+\frac{d}{p}$, it follows that
$$2^{js}\|R^{3}_{j}\|_{L^{p}_{x}(\mathcal{L}^{p})}\leq C c_{j}\|\nabla u \|_{B^{\frac{d}{p}}_{p,\infty}}\| \psi\|_{B^{s}_{p,r}(\mathcal{L}^{p})}.$$
Bounds for $2^{js}\|R^{4}_{j}\|_{L^{p}_{x}(\mathcal{L}^{p})}$, $2^{js}\|R^{5}_{j}\|_{L^{p}_{x}(\mathcal{L}^{p})}$, $2^{js}\|R^{6}_{j}\|_{L^{p}_{x}(\mathcal{L}^{p})}$ and $2^{js}\|R^{7}_{j}\|_{L^{p}_{x}(\mathcal{L}^{p})}$ are similar. We only treat with $R^{4}_{j}$.
Defining $\widetilde{\Delta}_{j'}=\Delta_{j'-1}+\Delta_{j'}+\Delta_{j'+1},$  we have
$$R^{4}_{j}=\sum_{|j-j'|\leq2}\partial_{k}(\Delta_{j}\widetilde{\Delta}_{j'}\psi\Delta_{j'}\widetilde{u}^{k}).$$
Hence, using (3.5), (3.6) and H\"{o}lder's inequality, we deduce that
$$2^{js}\|R^{4}_{j}\|_{L^{p}_{x}(\mathcal{L}^{p})}\leq C c_{j}\|\nabla u  \|_{L^{\infty}}\|\psi\|_{B^{s}_{p,r}(\mathcal{L}^{p})}.$$
\textit{Bounds for $2^{js}\|R^{8}_{j}\|_{L^{p}_{x}(\mathcal{L}^{p})}$.} A direct calculation yields
$$R^{8}_{j}=\sum_{|j-j'|\leq1}[\Delta_{j},\Delta_{-1}u]\cdot\nabla\Delta_{j'}\psi.$$
So by Lemma 3.16 we have
$$2^{js}\|R^{8}_{j}\|_{L^{p}_{x}(\mathcal{L}^{p})}\leq C\sum_{|j-j'|\leq 1}\|\nabla \Delta_{-1}u \|_{L^{\infty}}2^{j's}\|\Delta_{j'}\psi\|_{L^{p}_{x}(\mathcal{L}^{p})}\leq C c_{j}\|\nabla u  \|_{L^{\infty}}\|\psi\|_{B^{s}_{p,r}(\mathcal{L}^{p})}.$$
So combining the bounds for $R^{1}_{j}$ to $R^{8}_{j}$ yields the result.
\end{proof}
\end{lemm}

\begin{rema}
For homogeneous Besov spaces, Propositions 3.6, 3.7, 3.8, 3.15 and Lemmas 3.16, 3.18 still hold true. The proofs are similar, just replace $\Delta_j$ by $\dot\Delta_j$.
\end{rema}

\section{Linear problem and a priori estimates}
In this section we will consider the following linearized equations for (1.2):
\begin{align}
\left\{
\begin{array}{ll}
\partial_{t}u+(v\cdot\nabla)u+-\nu\Delta u+\nabla{P}=f, ~~ div u=0,\\[1ex]
\partial_{t}\psi+(v\cdot\nabla)\psi=div_{R}[-\nabla{v}R\psi+\psi_{\infty}\nabla_{R}\displaystyle\frac{\psi}{\psi_{\infty}}],  \\[1ex]
\tau_{ij}=\int_{B}R_{i}\otimes\nabla_{j}\mathcal{U}\psi dR, \\[1ex]
u|_{t=0}=u_{0}, \psi|_{t=0}=\psi_{0} ,\\[1ex]
\psi_{\infty}\nabla_{R}\displaystyle\frac{\psi}{\psi_{\infty}}\cdot n=0 ~~~~ on ~~~~ \partial B(0,R_{0}) .\\[1ex]
\end{array}
\right.
\end{align}

\subsection{Solutions to the linear equations in R}
Using Proposition 3.9 proved by Masmoudi in $[5]$, we can solve the following linear problem in $R$.
\begin{prop}
Assume that $A(t)\in C([0,T])$ is a matrix-valued function and $\psi_{0}\in\mathcal{L}^{p}$ with $p\in[2,+\infty)$, then
\begin{align}
\left\{
\begin{array}{ll}
\partial_{t}\psi=div_{R}[-A(t)R\psi+\psi_{\infty}\nabla_{R}\displaystyle\frac{\psi}{\psi_{\infty}}],  \\[1ex]
\psi|_{t=0}=\psi_{0} ,\\[1ex]
\psi_{\infty}\nabla_{R}\displaystyle\frac{\psi}{\psi_{\infty}}\cdot n=0 ~~~~ on ~~~~ \partial B(0,R_{0}) ,\\[1ex]
\end{array}
\right.
\end{align}
has a unique weak solution $\psi$ in $C([0,T];\mathcal{L}^{p})$. Moreover, we have
$$\sup_{t\in[0,T]}\int_{B}\bigg|\frac{\psi}{\psi_{\infty}}\bigg|^{p}\psi_{\infty}dR+\frac{2(p-1)}{p}\int^{T}_{0}\int_{B}\bigg|\nabla_{R}\bigg(\frac{\psi}{\psi_{\infty}}\bigg)^{\frac{p}{2}}\bigg|^{2}\psi_{\infty}dRdt\leq Ce^{C\int^{t}_{0}A^{2}(t')dt'}\int_{B}\bigg|\frac{\psi_{0}}{\psi_{\infty}}\bigg|^{p}\psi_{\infty}dR.$$
\begin{proof}
Firstly we smooth out the initial data $\psi_{0}$. Since $C^{\infty}(B)$ is dense in $\mathcal{L}^{p}$, it follows that there exists $\psi^{N}_{0}\in C^{\infty}(B)$  such that $$\psi^{N}_{0}\rightarrow\psi_{0}~~in~~ \mathcal{L}^{p}.$$
Assume that $\psi^{N}$ is the solution of
\begin{align}
\left\{
\begin{array}{ll}
\partial_{t}\psi^{N}=div_{R}[-A(t)R\psi^{N}+\psi_{\infty}\nabla_{R}\frac{\psi^{N}}{\psi_{\infty}}],  \\[1ex]
\psi^{N}|_{t=0}=\psi^{N}_{0} ,\\[1ex]
\psi_{\infty}\nabla_{R}\frac{\psi^{N}}{\psi_{\infty}}\cdot n=0 ~~~~ on ~~~~ \partial B(0,R_{0}) .\\[1ex]
\end{array}
\right.
\end{align}
As $p\geq2$, we have  $\mathcal{L}^{p}\hookrightarrow\mathcal{L}^{2}$. So by Proposition 3.9 in [1], we can exactly find $\psi^{N}\in C^{\infty}$. By multiplying both sides of (4.3) by $sgn(\psi^{N})\bigg|\displaystyle\frac{\psi^{N}}{\psi^{p-1}_{\infty}}\bigg|^{p-1}$ and integrating over $B$, we deduce that
\begin{align*}
\frac{1}{p}\partial_{t}\int_{B}\bigg|\frac{\psi^{N}}{\psi_{\infty}}\bigg|^{p}\psi_{\infty}dR&=(p-1)\bigg[\int_{B}A(t)R\bigg|\frac{\psi^{N}}{\psi_{\infty}}\bigg|^{p-1}\nabla_{\small R}\frac{\psi^{N}}{\psi_{\infty}}\psi_{\infty}dR-\int_{B}\bigg|\frac{\psi^{N}}{\psi_{\infty}}\bigg|^{p-2}\bigg(\nabla_{\small R}\frac{\psi^{N}}{\psi_{\infty}}\bigg)^{2}\psi_{\infty}dR\bigg] \\
&\leq (p-1)\bigg[C|A(t)|^{2}\int_{B}\bigg|\frac{\psi^{N}}{\psi_{\infty}}\bigg|^{p}\psi_{\infty}dR-\frac{1}{2}\int_{B}\bigg|\frac{\psi^{N}}{\psi_{\infty}}\bigg|^{p-2}\bigg(\nabla_{\small R}\frac{\psi^{N}}{\psi_{\infty}}\bigg)^{2}\psi_{\infty}dR\bigg] \\
&=C(p-1)|A(t)|^{2}\int_{B}\bigg|\frac{\psi^{N}}{\psi_{\infty}}\bigg|^{p}\psi_{\infty}dR-\frac{2(p-1)}{p^{2}}\int_{B}\bigg|\nabla_{R}\bigg(\frac{\psi^{N}}{\psi_{\infty}}\bigg)^{\frac{p}{2}}\bigg|^{2}\psi_{\infty}dR.
\end{align*}
Then we obtain
\begin{multline*}
\int_{B}\bigg|\frac{\psi^{N}}{\psi_{\infty}}\bigg|^{p}\psi_{\infty}dR+
\frac{2(p-1)}{p}\int^{t}_{0}\int_{B}\bigg|\nabla_{R}\bigg(\frac{\psi^{N}}{\psi_{\infty}}\bigg)^{\frac{p}{2}}\bigg|^{2}\psi_{\infty}dR \\
\leq\int_{B}\bigg|\frac{\psi^{N}_{0}}{\psi_{\infty}}\bigg|^{p}\psi_{\infty}dR
+Cp(p-1)\int^{t}_{0}\int_{B}|A(t')|^{2}\bigg|\frac{\psi^{N}}{\psi_{\infty}}\bigg|^{p}\psi_{\infty}dRdt.
\end{multline*}
Using Gronwall's inequality, we have
\begin{align}
\int_{B}\bigg|\frac{\psi^{N}}{\psi_{\infty}}\bigg|^{p}\psi_{\infty}dR+
\frac{2(p-1)}{p}\int^{t}_{0}\int_{B}\bigg|\nabla_{R}\bigg(\frac{\psi^{N}}{\psi_{\infty}}\bigg)^{\frac{p}{2}}\bigg|^{2}\psi_{\infty}dR \leq C e^{C\int^{t}_{0}A^{2}(t')dt'}\int_{B}\bigg|\frac{\psi^{N}_{0}}{\psi_{\infty}}\bigg|^{p}\psi_{\infty}dR.
\end{align}
Because (4.3) is a linear equation, $\psi^{N}-\psi^{M}$ has the same estimate as above, that is
$$\sup_{t\in [0,T]}\|\psi^{N}-\psi^{M}\|_{\mathcal{L}^{p}}\leq C e^{C\int^{T}_{0}A^{2}(t')dt'}\|\psi^{N}_{0}-\psi^{M}_{0}\|_{\mathcal{L}^{p}}.$$
Then $\psi^{N}$ is a Cauchy sequence in $L^{\infty}_{T}\mathcal{L}^{p}$.
There exists a $\psi\in L^{\infty}([0,T];\mathcal{L}^{p})$ such that
$$\psi^{N}\rightarrow\psi~~~in~~~L^{\infty}([0,T];\mathcal{L}^{p}).$$
Passing to the limit in the equation (4.3) we can see that $\psi$ is a solution of (4.2). Passing to the limit in (4.4), we obtain
$$\sup_{t\in[0,T]}\int_{B}\bigg|\frac{\psi}{\psi_{\infty}}\bigg|^{p}\psi_{\infty}dR+\frac{2(p-1)}{p}\int^{T}_{0}\int_{B}\bigg|\nabla_{R}\bigg(\frac{\psi}{\psi_{\infty}}\bigg)^{\frac{p}{2}}\bigg|^{2}\psi_{\infty}dRdt\leq Ce^{C\int^{t}_{0}A^{2}(t')dt'}\int_{B}\bigg|\frac{\psi_{0}}{\psi_{\infty}}\bigg|^{p}\psi_{\infty}dR.$$
Assume that $\psi,\phi$ are two solutions of (4.2) with the same initial data. From the above estimate we have
$$\|\psi-\phi\|_{\mathcal{L}^{p}}\leq 0,$$
which leads to the uniqueness.\\
Finally we shall prove that $\psi\in C([0,T];\mathcal{L}^{p})$. For any $t,s\in[0,T]$
\begin{align*}
\|\psi(t)-\psi(s)\|_{\mathcal{L}^{p}} &\leq \int^{t}_{s}\frac{1}{p}\bigg|\partial_{t}\int_{B}\bigg|\frac{\psi}{\psi_{\infty}}\bigg|^{p}\psi_{\infty}\bigg|dRdt' \\
&\leq(p-1)\bigg|C\int^{t}_{s}\int_{B}|A(t')|^{2}\bigg|\frac{\psi}{\psi_{\infty}}\bigg|^{p}\psi_{\infty}dR+\frac{1}{2}\int^{t}_{s}\int_{B}\bigg|\frac{\psi}{\psi_{\infty}}\bigg|^{p-2}\bigg(\nabla_{\small R}\frac{\psi}{\psi_{\infty}}\bigg)^{2}\psi_{\infty}dRdt'\bigg| \\
&\leq C(|t-s|\sup_{t\in[0,T]}A^{2}(t)\|\psi\|_{{\mathcal{L}^{p}}}+\int^{t}_{s}\int_{B}\bigg|\nabla_{R}\bigg(\frac{\psi}{\psi_{\infty}}\bigg)^{\frac{p}{2}}\bigg|^{2}\psi_{\infty}dRdt').
\end{align*}
Since $$\int^{T}_{0}\int_{B}\bigg|\nabla_{R}\bigg(\frac{\psi}{\psi_{\infty}}\bigg)^{\frac{p}{2}}\bigg|^{2}\psi_{\infty}dRdt<+\infty,$$
it then follows that
$$\|\psi(t)-\psi(s)\|_{\mathcal{L}^{p}}\rightarrow 0,~as~t\rightarrow s.$$
Then $\psi\in C([0,T];\mathcal{L}^{p})$.
\end{proof}
\end{prop}

Now we give a priori estimate for the Fokker-Planck equation.
\begin{lemm}
Assume that $\psi_{0}\in B^{\sigma}_{p,r}(\mathcal{L}^{p})$ and $f,g\in L^{2}([0,T];B^{\sigma}_{p,r}(\mathcal{L}^{p}))$, $u\in L^{\infty}([0,T];B^{s}_{p,r})\cap L^{2}([0,T];B^{s+1}_{p,r})$, and  $div~u =0$, where $s>1+\frac{d}{p},~s-1\leq \sigma\leq s,~p\in[2,+\infty),~r\geq p.$ If $\psi$  is a solution of
\begin{align}
\left\{
\begin{array}{ll}
\partial_{t}\psi+(u\cdot\nabla)\psi=div_{R}[-\nabla{u}R\psi+\psi_{\infty}\nabla_{R}\displaystyle\frac{\psi}{\psi_{\infty}}+g]+f,  \\[1ex]
\psi|_{t=0}=\psi_{0} ,\\[1ex]
\psi_{\infty}\nabla_{R}\displaystyle\frac{\psi}{\psi_{\infty}}\cdot n=0 ~~~~ on ~~~~ \partial B(0,R_{0}) .\\[1ex]
\end{array}
\right.
\end{align}
Then we have the following estimate
$$\sup_{t\in[0,T]}\|\psi\|_{B^{\sigma}_{p,r}(\mathcal{L}^{p})}+\|\psi\|_{E^{\sigma}_{p,r}(T)}\leq Ce^{CU(T)}\bigg(\|\psi_{0}\|_{B^{\sigma}_{p,r}(\mathcal{L}^{p})}+\bigg(\int^{T}_{0}e^{-CU(t')}(\|g(t')\|^{2}_{B^{\sigma}_{p,r}(\mathcal{L}^{p})}+\|f(t')\|^{2}_{B^{\sigma}_{p,r}(\mathcal{L}^{p})}dt'\bigg)^{\frac{1}{2}}\bigg),$$
where $U(t)=\int^{t}_{0}(\|\nabla u\|^{2}_{B^{s}_{p,r}}+1)dt'$.
\begin{proof}
Applying $\Delta_{j}$ to (4.5) yields
\begin{align}
\left\{
\begin{array}{ll}
\partial_{t}\Delta_{j}\psi+(u\cdot\nabla)\Delta_{j}\psi=div_{R}[-\Delta_{j}(\nabla{u}R\psi)+\psi_{\infty}\nabla_{R}\Delta_{j}\displaystyle\frac{\psi}{\psi_{\infty}}+\Delta_{j}g]+\Delta_{j}f+R_{j},  \\[1ex]
\Delta_{j}\psi|_{t=0}=\Delta_{j}\psi_{0} ,\\[1ex]
\psi_{\infty}\nabla_{R}\Delta_{j}\displaystyle\frac{\psi}{\psi_{\infty}}\cdot n=0 ~~~~ on ~~~~ \partial B(0,R_{0}) ,\\[1ex]
\end{array}
\right.
\end{align}
where $R_{j}=[u\cdot\nabla,\Delta_{j}]\psi.$  Since $u\in C([0,T];C^{0,1})$, we can define the flow of $u$, namely $\Phi(t,x)$  such that
\begin{align*}
\left\{
\begin{array}{ll}
\partial_{t}\Phi(t,x)=u(t,\Phi(t,x)),  \\[1ex]
\Phi(t,x)|_{t=0}=x .\\[1ex]
\end{array}
\right.
\end{align*}
For any function $a(t,x,R)$, we let $\widetilde{a}(t,x,R)=a(t,\Phi(t,x),R)$. Then (4.6) is equivalent to
\begin{align}
\left\{
\begin{array}{ll}
\partial_{t}\widetilde{\Delta_{j}\psi}=div_{R}[-\widetilde{\Delta_{j}(\nabla{u}R\psi)}+\psi_{\infty}\nabla_{R}\widetilde{\Delta_{j}\displaystyle\frac{\psi}{\psi_{\infty}}}+\widetilde{\Delta_{j}g}]+\widetilde{\Delta_{j}f}+\widetilde{R_{j}},  \\[1ex]
\widetilde{\Delta_{j}\psi}|_{t=0}=\Delta_{j}\psi_{0} ,\\[1ex]
\psi_{\infty}\nabla_{R}\widetilde{\Delta_{j}\displaystyle\frac{\psi}{\psi_{\infty}}}\cdot n=0 ~~~~ on ~~~~ \partial B(0,R_{0}) .\\[1ex]
\end{array}
\right.
\end{align}
By multiplying both sides of (4.7) by $sgn(\widetilde{\Delta_{j}\psi})\displaystyle\frac{|\widetilde{\Delta_{j}\psi}|^{p-1}}{\psi^{p-1}_{\infty}}$ and integrating over $\mathbb{R}^{d}\times B$, we have
\begin{align*}
\frac{1}{p}\partial_{t}\int_{\mathbb{R}^{d}\times B}\bigg|\frac{\widetilde{\Delta_{j}\psi}}{\psi_{\infty}}\bigg|^{p}\psi_{\infty}dxdR&=(p-1)\bigg[\int_{\mathbb{R}^{d}\times B}\bigg[\widetilde{\Delta_{j}\bigg(\nabla{u}R\frac{\psi}{\psi_{\infty}}\bigg)}+\widetilde{\Delta_{j}g}\bigg]\bigg|\frac{\widetilde{\Delta_{j}\psi}}{\psi_{\infty}}\bigg|^{p-2}\nabla_{\small R}\frac{\widetilde{\Delta_{j}\psi}}{\psi_{\infty}}\psi_{\infty}dxdR \\
&-\int_{\mathbb{R}^{d}\times B}\bigg|\frac{\widetilde{\Delta_{j}\psi}}{\psi_{\infty}}\bigg|^{p-2}\bigg(\nabla_{\small R}\frac{\widetilde{\Delta_{j}\psi}}{\psi_{\infty}}\bigg)^{2}\psi_{\infty}dxdR\bigg]
+\int_{\mathbb{R}^{d}\times B}(\widetilde{\Delta_{j}f}+\widetilde{R_{j}})\bigg|\frac{\widetilde{\Delta_{j}\psi}}{\psi_{\infty}}\bigg|^{p-1}dxdR\\
&\leq (p-1)\bigg[C\int_{\mathbb{R}^{d}\times B}\bigg[\bigg(\widetilde{\Delta_{j}\bigg(\nabla{u}R\frac{\psi}{\psi_{\infty}}\bigg)}\bigg)^{2}+|\widetilde{\Delta_{j}g}|^{2}\bigg]\bigg|\frac{\widetilde{\Delta_{j}\psi}}{\psi_{\infty}}\bigg|^{p-2}\psi_{\infty}dxdR\\
&-\frac{1}{2}\int_{\mathbb{R}^{d}\times B}\bigg|\frac{\widetilde{\Delta_{j}\psi}}{\psi_{\infty}}\bigg|^{p-2}\bigg(\nabla_{\small R}\frac{\widetilde{\Delta_{j}\psi}}{\psi_{\infty}}\bigg)^{2}\psi_{\infty}dxdR\bigg]
+\int_{\mathbb{R}^{d}\times B}(\widetilde{\Delta_{j}f}+\widetilde{R_{j}})\bigg|\frac{\widetilde{\Delta_{j}\psi}}{\psi_{\infty}}\bigg|^{p-1}dxdR\\
&=C(p-1)\int_{\mathbb{R}^{d}\times B}\bigg[\bigg(\widetilde{\Delta_{j}\bigg(\nabla{u}R\frac{\psi}{\psi_{\infty}}\bigg)}\bigg)^{2}+|\widetilde{\Delta_{j}g}|^{2}\bigg]\bigg|\frac{\widetilde{\Delta_{j}\psi}}{\psi_{\infty}}\bigg|^{p-2}\psi_{\infty}dxdR\\
&-\frac{2(p-1)}{p^{2}}\int_{\mathbb{R}^{d}\times B}\bigg|\nabla_{R}\bigg(\frac{\widetilde{\Delta_{j}\psi}}{\psi_{\infty}}\bigg)^{\frac{p}{2}}\bigg|^{2}\psi_{\infty}dxdR
+\int_{\mathbb{R}^{d}\times B}(\widetilde{\Delta_{j}f}+\widetilde{R_{j}})\bigg|\frac{\widetilde{\Delta_{j}\psi}}{\psi_{\infty}}\bigg|^{p-1}dxdR.
\end{align*}
Since $u$ is divergence-free, then the flow of $u$ is measure-preserving. From the above inequality we obtain
\begin{multline}
\frac{1}{p}\partial_{t}\int_{\mathbb{R}^{d}\times B}\bigg|\frac{\Delta_{j}\psi}{\psi_{\infty}}\bigg|^{p}\psi_{\infty}dxdR+\frac{2(p-1)}{p^{2}}\int_{\mathbb{R}^{d}\times B}\bigg|\nabla_{R}\bigg(\frac{\Delta_{j}\psi}{\psi_{\infty}}\bigg)^{\frac{p}{2}}\bigg|^{2}\psi_{\infty}dxdR \leq \\ C\int_{\mathbb{R}^{d}\times B}\bigg[\bigg(\Delta_{j}\bigg(\nabla{u}R\frac{\psi}{\psi_{\infty}}\bigg)\bigg)^{2}+|\Delta_{j}g|^{2}\bigg]\bigg|\frac{\Delta_{j}\psi}{\psi_{\infty}}\bigg|^{p-2}\psi_{\infty}dxdR
+\int_{\mathbb{R}^{d}\times B}(\Delta_{j}f+R_{j})\bigg|\frac{\Delta_{j}\psi}{\psi_{\infty}}\bigg|^{p-1}dxdR.
\end{multline}
Using H\"{o}lder's inequality, we deduce that
\begin{align*}
\partial_{t}\|\Delta_{j}\psi\|^{p}_{L^{p}_{x}(\mathcal{L}^{p})} &\leq C [(\|\Delta_{j}(\nabla{u}R\psi)\|^{2}_{L^{p}_{x}(\mathcal{L}^{p})}+\|\Delta_{j}g\|^{2}_{L^{p}_{x}(\mathcal{L}^{p})})\|\Delta_{j}\psi\|^{p-2}_{L^{p}_{x}(\mathcal{L}^{p})}\\
&+(\|\Delta_{j}f\|_{L^{p}_{x}(\mathcal{L}^{p})}+\|R_{j}\|_{L^{p}_{x}(\mathcal{L}^{p})})\|\Delta_{j}\psi\|^{p-1}_{L^{p}_{x}(\mathcal{L}^{p})}].
\end{align*}
Thus we have
\begin{align*}
\partial_{t}\|\Delta_{j}\psi\|^{2}_{L^{p}_{x}(\mathcal{L}^{p})} \leq C [\|\Delta_{j}(\nabla{u}R\psi)\|^{2}_{L^{p}_{x}(\mathcal{L}^{p})}+\|\Delta_{j}g\|^{2}_{L^{p}_{x}(\mathcal{L}^{p})})+(\|\Delta_{j}f\|_{L^{p}_{x}(\mathcal{L}^{p})}+\|R_{j}\|_{L^{p}_{x}(\mathcal{L}^{p})})\|\Delta_{j}\psi\|_{L^{p}_{x}(\mathcal{L}^{p})}].
\end{align*}
Integrating over $[0,t]$, we deduce that
\begin{multline}
\|\Delta_{j}\psi\|^{2}_{L^{p}_{x}(\mathcal{L}^{p})} \leq \|\Delta_{j}\psi_{0}\|^{2}_{L^{p}_{x}(\mathcal{L}^{p})}+
\int^{t}_{0}\|\Delta_{j}(\nabla{u}R\psi)\|^{2}_{L^{p}_{x}(\mathcal{L}^{p})}+\|\Delta_{j}g\|^{2}_{L^{p}_{x}(\mathcal{L}^{p})}\\
+(\|\Delta_{j}f\|_{L^{p}_{x}(\mathcal{L}^{p})}+
\|R_{j}\|_{L^{p}_{x}(\mathcal{L}^{p})})\|\Delta_{j}\psi\|_{L^{p}_{x}(\mathcal{L}^{p})}dt'.
\end{multline}
Let $(c_{j})_{j\geq-1}$ denotes an element of the unit sphere of $l^{r}$. By Lemma (3.17) we have
\begin{align}
\|R_{j}\|_{L^{p}_{x}(\mathcal{L}^{p})} \leq C c_{j}2^{-j\sigma} \|\nabla u\|_{B^{s-1}_{p,r}}\|\psi\|_{B^{\sigma}_{p,r}(\mathcal{L}^{p})}.
\end{align}
Plugging (4.10) into (4.9) and by the Cauchy-Schwarz inequality we have
\begin{align}
\nonumber \|\Delta_{j}\psi\|^{2}_{L^{p}_{x}(\mathcal{L}^{p})} &\leq \|\Delta_{j}\psi_{0}\|^{2}_{L^{p}_{x}(\mathcal{L}^{p})}+
\int^{t}_{0}\|\Delta_{j}(\nabla{u}R\psi)\|^{2}_{L^{p}_{x}(\mathcal{L}^{p})}+\|\Delta_{j}g\|^{2}_{L^{p}_{x}(\mathcal{L}^{p})}\\
&+\|\Delta_{j}f\|^{2}_{L^{p}_{x}(\mathcal{L}^{p})}+\|\Delta_{j}\psi\|^{2}_{L^{p}_{x}(\mathcal{L}^{p})}+
Cc_{j}2^{-j\sigma} \|\nabla u\|_{B^{s-1}_{p,r}}\|\psi\|_{B^{\sigma}_{p,r}(\mathcal{L}^{p})}\|\Delta_{j}\psi\|_{L^{p}_{x}(\mathcal{L}^{p})}dt'.
\end{align}
By multiplying both sides $2^{2j\sigma}$ and taking $l^{\frac{r}{2}}$ norm (here we use the fact that $r\geq2$), we deduce that
\begin{align*}
 \|\psi\|^{2}_{B^{\sigma}_{p,r}(\mathcal{L}^{p})} &\leq C\{\|\psi_{0}\|^{2}_{B^{\sigma}_{p,r}(\mathcal{L}^{p})}+
\int^{t}_{0}\|\nabla{u}R\psi\|^{2}_{B^{\sigma}_{p,r}(\mathcal{L}^{p})}+\|g\|^{2}_{B^{\sigma}_{p,r}(\mathcal{L}^{p})}
+\|f\|^{2}_{B^{s}_{p,r}(\mathcal{L}^{p})}+
(1+ \|\nabla u\|_{B^{s-1}_{p,r}})\|\psi\|^{2}_{B^{\sigma}_{p,r}(\mathcal{L}^{p})}dt'\}  \\
 &\leq C \{\|\psi_{0}\|^{2}_{B^{\sigma}_{p,r}(\mathcal{L}^{p})}+\int^{t}_{0}\|\nabla{u}\|^{2}_{B^{s}_{p,r}(\mathcal{L}^{p})}\|\psi\|^{2}_{B^{\sigma}_{p,r}(\mathcal{L}^{p})}
+(1+ \|\nabla u\|_{B^{\sigma-1}_{p,r}})\|\psi\|^{2}_{B^{\sigma}_{p,r}(\mathcal{L}^{p})}dt'\\
&+\int^{t}_{0}\|f\|^{2}_{B^{\sigma}_{p,r}(\mathcal{L}^{p})}+\|g\|^{2}_{B^{\sigma}_{p,r}(\mathcal{L}^{p})}dt'\}
\\
 &\leq C\{\|\psi_{0}\|^{2}_{B^{\sigma}_{p,r}(\mathcal{L}^{p})}+\int^{t}_{0}
(1+ \|\nabla u\|_{B^{s}_{p,r}}+\|\nabla{u}\|^{2}_{B^{s}_{p,r}(\mathcal{L}^{p})})\|\psi\|^{2}_{B^{\sigma}_{p,r}(\mathcal{L}^{p})}dt'\\
&+\int^{t}_{0}\|f\|^{2}_{B^{\sigma}_{p,r}(\mathcal{L}^{p})}+\|g\|^{2}_{B^{\sigma}_{p,r}(\mathcal{L}^{p})}dt'\}\\
&\leq C\{\|\psi_{0}\|^{2}_{B^{\sigma}_{p,r}(\mathcal{L}^{p})}+\int^{t}_{0}
\|f\|^{2}_{B^{\sigma}_{p,r}(\mathcal{L}^{p})}+\|g\|^{2}_{B^{\sigma}_{p,r}(\mathcal{L}^{p})}+(1+\|\nabla{u}\|^{2}_{B^{s}_{p,r}(\mathcal{L}^{p})})\|\psi\|^{2}_{B^{\sigma}_{p,r}(\mathcal{L}^{p})}dt'\}.
\end{align*}
Using Gronwall's inequality, we deduce that
\begin{align}
\sup_{t\in[0,T]}\|\psi\|_{B^{\sigma}_{p,r}(\mathcal{L}^{p})}\leq Ce^{CU(T)}\bigg(\|\psi_{0}\|_{B^{\sigma}_{p,r}(\mathcal{L}^{p})}+\bigg(\int^{T}_{0}e^{-CU(t')}(\|f(t')\|^{2}_{B^{\sigma}_{p,r}(\mathcal{L}^{p})}+\|g(t')\|^{2}_{B^{\sigma}_{p,r}(\mathcal{L}^{p})})dt'\bigg)^{\frac{1}{2}}\bigg), \end{align}
where $U(t)=\int^{t}_{0}\|\nabla u\|^{2}_{B^{s}_{p,r}}+1dt'.$
Now integrating (4.8) with respect to $t$ over $[0,T]$, we obtain
\begin{multline}
\displaystyle\int^{T}_{0}\displaystyle\int_{\mathbb{R}^{d}\times B}\bigg|\nabla_{R}\bigg(\Delta_{j}\displaystyle\frac{\psi}{\psi_{\infty}}\bigg)^{\frac{p}{2}}\bigg|^{2}\psi_{\infty}dxdRdt \leq C \|\Delta_{j}\psi_{0}\|^{p}_{L^{p}_{x}(\mathcal{L}^{p})}\\
+C\int^{T}_{0}(\|\Delta_{j}(\nabla{u}R\psi)\|^{2}_{L^{p}_{x}(\mathcal{L}^{p})}+\|\Delta_{j}g\|^{2}_{L^{p}_{x}(\mathcal{L}^{p})})\|\Delta_{j}\psi\|^{p-2}_{L^{p}_{x}(\mathcal{L}^{p})}\\
+(\|\Delta_{j}f\|_{L^{p}_{x}(\mathcal{L}^{p})}+\|R_{j}\|_{L^{p}_{x}(\mathcal{L}^{p})})\|\Delta_{j}\psi\|^{p-1}_{L^{p}_{x}(\mathcal{L}^{p})})dt.
\end{multline}
 Multiplying both sides of (4.13) by $2^{pj\sigma}$ and taking $l^{\frac{r}{p}}$-norm (here we use the fact that $r\geq p$), and using the inequality (4.10) we deduce that
\begin{align}
\|\psi\|_{E^{\sigma}_{p,r}(T)} \leq Ce^{CU(T)}\bigg(\|\psi_{0}\|_{B^{\sigma}_{p,r}(\mathcal{L}^{p})}+\bigg(\int^{T}_{0}e^{-CU(t')}(\|f(t')\|^{2}_{B^{\sigma}_{p,r}(\mathcal{L}^{p})}+\|g(t')\|^{2}_{B^{\sigma}_{p,r}(\mathcal{L}^{p})})dt'\bigg)^{\frac{1}{2}}\bigg).
\end{align}
Combining (4.12) and (4.14), we thus complete the proof.
\end{proof}
\end{lemm}

\begin{prop}
Assume that $\psi_{0}\in B^{s}_{p,r}(\mathcal{L}^{p})$ and $u\in L^{\infty}([0,T];B^{s}_{p,r})\cap L^{2}([0,T];B^{s+1}_{p,r})$, where $s>1+\frac{d}{p},~p\in[2,+\infty),~r\geq p.$ Then
\begin{align}
\left\{
\begin{array}{ll}
\partial_{t}\psi+(u\cdot\nabla)\psi=div_{R}[-\nabla{u}R\psi+\psi_{\infty}\nabla_{R}\displaystyle\frac{\psi}{\psi_{\infty}}],  \\[1ex]
\psi|_{t=0}=\psi_{0} ,\\[1ex]
\psi_{\infty}\nabla_{R}\displaystyle\frac{\psi}{\psi_{\infty}}\cdot n=0 ~~~~ on ~~~~ \partial B(0,R_{0}) ,\\[1ex]
\end{array}
\right.
\end{align}
has a unique solution $\psi$ in $C([0,T];B^{s}_{p,r}(\mathcal{L}^{p}))$ if $r<\infty$ or  in $C_{w}([0,T];B^{s}_{p,r}(\mathcal{L}^{p}))$ if $r=\infty$. Where $C_{w}([0,T];X)$ denotes the weak continuous space over $[0,T]$ on the Banach space $X$.
Moreover, we have the following estimate
$$\sup_{t\in[0,T]}\|\psi\|_{B^{s}_{p,r}(\mathcal{L}^{p})}+\|\psi\|_{E^s_{p,r}(T)}\leq Ce^{CU(T)}\|\psi_{0}\|_{B^{s}_{p,r}(\mathcal{L}^{p})},$$
where $U(t)=\int^{t}_{0}(\|\nabla u\|^{2}_{B^{s}_{p,r}}+1)dt'$.
\begin{proof}
As in Lemma 4.2, defining the flow of $u$, and letting $\widetilde{a(t,x,R)}=a(t,\Phi(t,x),R)$
 so we can write (4.15) as
\begin{align}
\left\{
\begin{array}{ll}
\partial_{t}\widetilde{\psi}=div_{R}[-\widetilde{\nabla{u}R\psi}+\psi_{\infty}\nabla_{R}\displaystyle\frac{\widetilde{\psi}}{\psi_{\infty}}],  \\[1ex]
\widetilde{\psi}|_{t=0}=\psi_{0} ,\\[1ex]
\psi_{\infty}\nabla_{R}\frac{\widetilde{\psi}}{\psi_{\infty}}\cdot n=0 ~~~~ on ~~~~ \partial B(0,R_{0}) .\\[1ex]
\end{array}
\right.
\end{align}
Then using Proposition (4.1) for each fixed $x$, we deduce the existence and uniqueness of $\psi(t,x,R)\in C([0,T];\mathcal{L}^{p})$. Thus using Lemma 4.2 with $f,g=0,~\sigma=s$, we get the following estimate
$$\sup_{t\in[0,T]}\|\psi\|_{B^{s}_{p,r}(\mathcal{L}^{p})}+\|\psi\|_{E^{s}_{p,r}(T)}\leq Ce^{CU(T)}\|\psi_{0}\|_{B^{s}_{p,r}(\mathcal{L}^{p})}.$$
So we have $\psi\in L^{\infty}([0,T];B^{s}_{p,r}(\mathcal{L}^{p})$.
Since the equation  is a linear equation, using the above estimate we thus prove the uniqueness
in the space  $\psi\in L^{\infty}([0,T];B^{s}_{p,r}(\mathcal{L}^{p})$. \\
Next we will check that $\psi\in C([0,T];B^{s}_{p,r}(\mathcal{L}^{p})$ when $r<\infty$.
Applying $\Delta_{j}$ to (4.15) and by a similar calculation as in Lemma 4.2 with $f,g=0$, we can deduce that
\begin{align*}
\partial_{t}\|\Delta_{j}\psi\|^{2}_{L^{p}_{x}(\mathcal{L}^{p})} \leq C \big(\|\Delta_{j}(\nabla{u}R\psi)\|^{2}_{L^{p}_{x}(\mathcal{L}^{p})}+\|R_{j}\|_{L^{p}_{x}(\mathcal{L}^{p})}\|\Delta_{j}\psi\|_{L^{p}_{x}(\mathcal{L}^{p})}\big).
\end{align*}
Let $(c_{j})_{j\geq-1}$ denote an element of the unit sphere of $l^{r}$. By the definition of $B^{s}_{p,r}(\mathcal{L}^{p})$ and Lemma 3.17. We see that
\begin{align}
\nonumber \partial_{t}\|\Delta_{j}\psi\|^{2}_{L^{p}_{x}(\mathcal{L}^{p})} &\leq C c^{2}_{j}2^{-2js} \big(\|(\nabla{u}R\psi)\|^{2}_{B^{s}_{p,r}(\mathcal{L}^{p})}+\|\nabla{u}\|_{B^{s-1}_{p,r}}\|\psi\|_{B^{s}_{p,r}(\mathcal{L}^{p})}\big)\\
&\leq C c^{2}_{j}2^{-2js}\big[(\|\nabla{u}\|^{2}_{B^{s}_{p,r}}+\|\nabla{u}\|_{B^{s-1}_{p,r}})\|\psi\|_{B^{s}_{p,r}(\mathcal{L}^{p})}\big]. \end{align}
Due to $\psi\in L^{\infty}([0,T];B^{s}_{p,r}(\mathcal{L}^{p}))$, for any $\varepsilon>0$, there exists $N$ such that
\begin{align}
\sup_{t\in[0,T]}\sum_{j\geq N}2^{jsr}\|\Delta_{j}\psi\|^{r}_{L^{p}_{x}(\mathcal{L}^{p})} <\frac{\varepsilon}{4},
\end{align}
hence, for any $t_{1},t_{2}\in[0,T]$, we have
\begin{align*}
\|\psi(t_{1})-\psi(t_{2})\|_{B^{s}_{p,r}(\mathcal{L}^{p})} &\leq  (\sum_{-1\leq j< N}2^{jsr}\|\Delta_{j}\psi(t_{1})-\Delta_{j}\psi(t_{2})\|^{r}_{L^{p}_{x}(\mathcal{L}^{p})} +2\sup_{t\in[0,T]}\sum_{j\geq N}2^{jsr}\|\Delta_{j}\psi\|^{r}_{L^{p}_{x}(\mathcal{L}^{p})})^{\frac{1}{r}} \\
&\leq  (\sum_{-1\leq j< N}2^{jsr}\|\Delta_{j}\psi(t_{1})-\Delta_{j}\psi(t_{2})\|^{r}_{L^{p}_{x}(\mathcal{L}^{p})})^{\frac{1}{r}}+\frac{\varepsilon}{2}\\ &\leq (\sum_{-1\leq j< N}2^{jsr}\int^{t_{2}}_{t_{1}}\partial_{t}\|\Delta_{j}\psi(t')\|^{r}_{L^{p}_{x}(\mathcal{L}^{p})}dt')^{\frac{1}{r}}+\frac{\varepsilon}{2}\\
&\leq  (\sum_{-1\leq j< N}\frac{r}{2}2^{jsr}\int^{t_{2}}_{t_{1}}\partial_{t}\|\Delta_{j}\psi(t')\|^{2}_{L^{p}_{x}(\mathcal{L}^{p})}\|\Delta_{j}\psi(t')\|^{r-2}_{L^{p}_{x}(\mathcal{L}^{p})}dt')^{\frac{1}{r}}+\frac{\varepsilon}{2}.
\end{align*}
By the inequality (4.17) and the definition of $B^{s}_{p,r}(\mathcal{L}^{p})$, we obtain
\begin{align*}
\|\psi(t_{1})-\psi(t_{2})\|_{B^{s}_{p,r}(\mathcal{L}^{p})} &\leq  (\sum_{-1\leq j< N}\frac{r}{2}\int^{t_{2}}_{t_{1}}c^{r}_{j}(\|\nabla u\|^{2}_{B^{s}_{p,r}}+\|\nabla u\|_{B^{s-1}_{p,r}})dt')^{\frac{1}{r}}\sup_{t\in[0,T]}\|\psi(t)\|_{B^{s}_{p,r}(\mathcal{L}^{p})}+\frac{\varepsilon}{2}\\
&\leq C(N+1)(\int^{t_{2}}_{t_{1}}\|\nabla u\|^{2}_{B^{s}_{p,r}}+\|\nabla u\|_{B^{s-1}_{p,r}}dt')^{\frac{1}{r}}\sup_{t\in[0,T]}\|\psi(t)\|_{B^{s}_{p,r}(\mathcal{L}^{p})}+\frac{\varepsilon}{2}.
\end{align*}
Since $u\in L^{\infty}([0,T];B^{s}_{p,r})\cap L^{2}([0,T];B^{s+1}_{p,r})$, it follows that $\nabla u\in L^{\infty}([0,T];B^{s-1}_{p,r})\cap L^{2}([0,T];B^{s}_{p,r})$. Then we have
$$\|\psi(t_{1})-\psi(t_{2})\|_{B^{s}_{p,r}(\mathcal{L}^{p})}\rightarrow 0~as~t_1\rightarrow t_2.$$
Thus we prove that $\psi \in C([0,T];B^{s}_{p,r}(\mathcal{L}^{p}))$ when $r<\infty$.
Finally we consider $r=\infty$. Since $s>1+\frac{d}{p}$, there exist $s'<s$ such that $s'>1+\frac{d}{p}$.
Using the fact that $B^{s}_{p,\infty}\hookrightarrow B^{s'}_{p,r}$ for any $r\geq 1$. Fix some $r\in [p,\infty)$,  by the previous argument, we have  $\psi\in C([0,T];B^{s'}_{p,r}(\mathcal{L}^{p}))$. Now we need to check that $\psi\in C_{w}([0,T];B^{s'}_{p,r})$. Indeed , for any $\phi \in S(R^{d};C^{\infty}_{0}(B))$
$$\langle\psi(t_{1})-\psi(t_{2}),\phi\rangle\leq C\|\psi(t_{1})-\psi(t_{2})\|_{B^{s'}_{p,r}(\mathcal{L}^{p}))}\|\phi\|_{B^{-s'}_{p',r'}(L^{p'};\psi^{1-p}_{\infty} dR)}.$$
We thus complete the proof.
\end{proof}
\end{prop}

\subsection{A priori estimates for the Navier-Stokes equations}
In the following lemma, we give a priori estimates for the linear Navier-Stokes equations.
\begin{lemm}
Assume that $u_{0}\in B^{s}_{p,r}$, $v\in L^{\infty}([0,T];B^{s}_{p,r})$ with $div~v=0$, $f\in L^{2}([0,T];B^{s-1}_{p,r})$ and $\psi\in L^{\infty}([0,T];B^{s}_{p,r})\cap E^{s}_{p,r}(T)$ where $s>1+\frac{d}{p},p\in[2;\infty),r\geq p$, $u$ is the solution of
\begin{align}
\left\{
\begin{array}{ll}
\partial_{t}u+(v\cdot\nabla)u-\nu\Delta u+\nabla{P}=div \tau+f, ~ div u=0,\\[1ex]
\tau_{ij}=\int_{B}R_{i}\otimes\nabla_{j}\mathcal{U}\psi dR, \\[1ex]
u|_{t=0}=u_{0} .\\[1ex]
\end{array}
\right.
\end{align}
Then for any $\varepsilon>0$,  we have following estimates:
$$\sup_{t\in[0,T]}\|u\|^{2}_{B^{s}_{p,r}}\leq Ce^{CV(T)}\bigg[\|u_{0}\|^{2}_{B^{s}_{p,r}}+(\nu^{-1}+T)(\varepsilon\|\psi\|^{2}_{E^{s}_{p,r}(T)}+\int^{T}_{0}e^{-CV(t')}(\|f\|^{2}_{B^{s-1}_{p,r}}+\|\psi\|^{2}_{B^{s}_{p,r}(\mathcal{L}^{p})})dt')\bigg],$$
$$\nu\int^{T}_{0}\|u\|^{2}_{B^{s+1}_{p,r}}dt \leq C(1+\nu T)e^{cV(T)}\bigg[\|u_{0}\|^{2}_{B^{s}_{p,r}}+(\varepsilon\|\psi\|^{2}_{E^{s}_{p,r}(T)}+\int^{T}_{0}e^{-cV(t)}(\|\psi\|^{2}_{B^{s}_{p,r}(\mathcal{L}^{p})}+\|f\|^{2}_{B^{s-1}_{p,r}})dt)\bigg],$$
where $V(t)=\int^{t}_{0}\|v\|^{2}_{B^{s}_{p,r}}dt'.$
\begin{proof}
Applying $\Delta_{j}$ to (4.19) yields
\begin{align}
\left\{
\begin{array}{ll}
\partial_{t}\Delta_{j}u-\nu\Delta \Delta_{j}u+\nabla{\Delta_{j}P}=div \Delta_{j}\tau+\Delta_{j}f-\Delta_{j}(v\cdot\nabla u),  \\[1ex]
div \Delta_{j}u=0,     \\[1ex]
\Delta_{j}u|_{t=0}=\Delta_{j}u_{0}. \\[1ex]
\end{array}
\right.
\end{align}
So we can write that
\begin{align}
\Delta_{j}u=e^{\nu t\Delta}\Delta_{j}u_{0}+\int^{t}_{0}e^{\nu (t-t')\Delta}(div \Delta_{j}\tau+\Delta_{j}f-\Delta_{j}(v\cdot\nabla u)-\nabla{\Delta_{j}P})dt'.
\end{align}
If $j\geq 0$,  by Lemma 3.2 we have
\begin{multline}
\|\Delta_{j}u\|_{L^{p}}\leq C\bigg(e^{-c\nu t2^{2j}}\|\Delta_{j}u_{0}\|_{L^{p}}+\int^{t}_{0}e^{-c\nu (t-t')2^{2j}}(\|div \Delta_{j}\tau\|_{L^{p}}\\
+\|\Delta_{j}f\|_{L^{p}}+\|\Delta_{j}(v\cdot\nabla u)\|_{L^{p}}+\|\nabla{\Delta_{j}P}\|_{L^{p}}dt')\bigg).
\end{multline}
Firstly we deal with the pressure term $\|\nabla{\Delta_{j}P}\|_{L^{p}}$. Taking $div$ for (4.20), we deduce that
$$\Delta \Delta_{j}P=div((div \Delta_{j}\tau+\Delta_{j}f-\Delta_{j}(v\cdot\nabla u)).$$
Then we see that $\Delta_{j}P$ is a solution of an elliptic equation. Thus, we obtain
$$\nabla\Delta_{j}P=\nabla (\Delta)^{-1}div ((div \Delta_{j}\tau+\Delta_{j}f-\Delta_{j}(v\cdot\nabla u)).$$
Thanks to $\nabla (\Delta)^{-1}div$ is a Calderon-Zygmund operator and $p<\infty$, we have
 $$\|\nabla{\Delta_{j}P}\|_{L^{p}}\leq C(\|div \Delta_{j}\tau\|_{L^{p}}+\|\Delta_{j}f\|_{L^{p}}+\|\Delta_{j}(v\cdot\nabla u)\|_{L^{p}}).$$
 Plugging into (4.22), we deduce that
  \begin{multline}
\|\Delta_{j}u\|_{L^{p}}\leq C\bigg(e^{-c\nu t2^{2j}}\|\Delta_{j}u_{0}\|_{L^{p}}+\int^{T}_{0}e^{-c\nu (t-t')2^{2j}}(\|div \Delta_{j}\tau\|_{L^{p}}
+\|\Delta_{j}f\|_{L^{p}}+\|\Delta_{j}(v\cdot\nabla u)\|_{L^{p}}dt'\bigg).
\end{multline}
Notice that
\begin{multline}
\int^{t}_{0}e^{-c\nu (t-t')2^{2j}}(\|div \Delta_{j}\tau\|_{L^{p}}
+\|\Delta_{j}f\|_{L^{p}}+\|\Delta_{j}(v\cdot\nabla u)\|_{L^{p}}dt'\\=e^{-c\nu t2^{2j}}1_{[0,T]}*(\|div \Delta_{j}\tau(t)\|_{L^{p}}
+\|\Delta_{j}f(t)\|_{L^{p}}+\|\Delta_{j}(v(t)\cdot\nabla u(t))\|_{L^{p}})1_{[0,T]}$$.
\end{multline}
Using Young's inequalty, we deduce that
\begin{align}
\nonumber \|\Delta_{j}u\|_{L^{p}}&\leq C\bigg[e^{-c\nu t2^{2j}}\|\Delta_{j}u_{0}\|_{L^{p}}+\bigg(\int^{T}_{0}e^{-c\nu t2^{2j}}dt\bigg)^{\frac{1}{2}}\bigg(\int^{T}_{0}(\|div \Delta_{j}\tau\|^{2}_{L^{p}}
+\|\Delta_{j}f\|^{2}_{L^{p}}+\|\Delta_{j}(v\cdot\nabla u)\|^{2}_{L^{p}}dt')\bigg)^{\frac{1}{2}}\bigg]\\
&\leq C \bigg[e^{-c\nu t2^{2j}}\|\Delta_{j}u_{0}\|_{L^{p}}+\nu^{-\frac{1}{2}}2^{-j}\bigg(\int^{T}_{0}(\|div \Delta_{j}\tau\|^{2}_{L^{p}}
+\|\Delta_{j}f\|^{2}_{L^{p}}+\|\Delta_{j}(v\cdot\nabla u)\|^{2}_{L^{p}}dt')\bigg)^{\frac{1}{2}}\bigg].
\end{align}
Since $div~v=0$, it follows that $v\cdot\nabla u=div (v\otimes u)$. And by Bernstein's inequality, we get
\begin{align}
\|\Delta_{j}u\|_{L^{p}}
&\leq C \bigg[e^{-c\nu t2^{2j}}\|\Delta_{j}u_{0}\|_{L^{p}}+\nu^{-\frac{1}{2}}\bigg(\int^{T}_{0}(\| \Delta_{j}\tau\|^{2}_{L^{p}}
+2^{-2j}\|\Delta_{j}f\|^{2}_{L^{p}}+\|\Delta_{j}(v\otimes u)\|^{2}_{L^{p}}dt')\bigg)^{\frac{1}{2}}\bigg].
\end{align}
That is
 \begin{align}
\|\Delta_{j}u\|^{2}_{L^{p}}
&\leq C \bigg[e^{-c\nu t2^{2j}}\|\Delta_{j}u_{0}\|^{2}_{L^{p}}+\nu^{-1}\bigg(\int^{T}_{0}(\| \Delta_{j}\tau\|^{2}_{L^{p}}
+2^{-2j}\|\Delta_{j}f\|^{2}_{L^{p}}+\|\Delta_{j}(v\otimes u)\|^{2}_{L^{p}}dt')\bigg)\bigg].
\end{align}
Now we deal with the stress tensor $\tau$. By Corollary (3.13), we obtain
$$\| \Delta_{j}\tau\|_{L^{p}}\leq C\bigg[\bigg(\varepsilon\displaystyle\int_{\mathbb{R}^{d}\times B}\bigg|\nabla_{R}\bigg(\Delta_{j}\displaystyle\frac{\psi}{\psi_{\infty}}\bigg)^{\frac{p}{2}}\bigg|^{2}\psi_{\infty}dxdR\bigg)^{\frac{1}{p}}+\|\Delta_{j}\psi\|_{L^{p}_{x}(\mathcal{L}^{p})}\bigg].$$
Plugging into (4.27), we deduce that
\begin{multline}
\|\Delta_{j}u\|^{2}_{L^{p}}
\leq C \bigg[e^{-c\nu t2^{2j}}\|\Delta_{j}u_{0}\|^{2}_{L^{p}}+\nu^{-1}\bigg(\int^{T}_{0}\bigg(\varepsilon\displaystyle\int_{\mathbb{R}^{d}\times B}\bigg|\nabla_{R}\bigg(\Delta_{j}\displaystyle\frac{\psi}{\psi_{\infty}}\bigg)^{\frac{p}{2}}\bigg|^{2}\psi_{\infty}dxdR\bigg)^{\frac{2}{p}}\\
+\|\Delta_{j}\psi\|^{2}_{L^{p}_{x}(\mathcal{L}^{p})}
+2^{-2j}\|\Delta_{j}f\|^{2}_{L^{p}}+\|\Delta_{j}(v\otimes u)\|^{2}_{L^{p}}dt')\bigg)\bigg].
\end{multline}
If $j=-1$, applying $\Delta_{-1}$ to (4.19) yields
\begin{align}
\left\{
\begin{array}{ll}
\partial_{t}\Delta_{-1}u+v\cdot\nabla \Delta_{-1}u-\nu\Delta \Delta_{-1}u+\nabla{\Delta_{-1}P}=div \Delta_{-1}\tau+\Delta_{-1}f-[v\cdot\nabla,\Delta_{-1}]u,  \\[1ex]
div \Delta_{-1}u=0,     \\[1ex]
\Delta_{-1}u|_{t=0}=\Delta_{-1}u_{0}. \\[1ex]
\end{array}
\right.
\end{align}
Multiplying both sides of (4.29) by $sgn(\Delta_{-1}u)|\Delta_{-1}u|^{p-1}$, and using the fact that $div~v=0$, we obtain
\begin{align}
\|\Delta_{-1}u\|_{L^{p}}&\leq \|\Delta_{-1}u_{0}\|_{L^{p}}+C\int^{T}_{0}\|div(\Delta_{-1})\tau\|_{L^{p}}+\|\Delta_{-1}f\|_{L^{p}}+\|[v\cdot\nabla,\Delta_{-1}]u\|_{L^{p}}dt'\\
\nonumber&\leq \|\Delta_{-1}u_{0}\|_{L^{p}}+
CT^{\frac{1}{2}}\bigg[\int^{T}_{0}\bigg(\varepsilon\displaystyle\int_{\mathbb{R}^{d}\times B}\bigg|\nabla_{R}\bigg(\Delta_{-1}\displaystyle\frac{\psi}{\psi_{\infty}}\bigg)^{\frac{p}{2}}\bigg|^{2}\psi_{\infty}dxdR\bigg)^{\frac{2}{p}}\\
\nonumber&+\|\Delta_{-1}\psi\|^{2}_{L^{p}_{x}(\mathcal{L}^{p})}+\|\Delta_{-1}f\|^{2}_{L^{p}}+\|\nabla v\|^{2}_{L^{\infty}}\|\Delta_{-1}u\|^{2}_{L^{p}}\bigg]^{\frac{1}{2}}.
\end{align}
Mutiplying both sidesof (4.28) by $2^{2js}$, and taking $l^{\frac{r}{2}}$-norm, and combining the inequality (4.30), we deduce that
\begin{align}
\|u\|^{2}_{B^{s}_{p,r}} \leq C\|u_{0}\|^{2}_{B^{s}_{p,r}}+(\nu^{-1}+T)\int^{T}_{0}\|\psi\|^{2}_{B^{s}_{p,r}(\mathcal{L}^{p})}+\|f\|^{2}_{B^{s-1}_{p,r}}+\|v\otimes u\|^{2}_{B^{s}_{p,r}}+\|\nabla v\|^{2}_{L^{\infty}}\|u\|^{2}_{B^{s}_{p,r}}dt+\varepsilon\|\psi\|_{E^{s}_{p,r}(T)}.
\end{align}
Since $s>1+\frac{d}{p}$, it follows that $B^{s}_{p,r}\hookrightarrow C^{0,1}$ is an algebra. Then we have
\begin{align}
\sup_{t\in[0,T]}\|u\|^{2}_{B^{s}_{p,r}} \leq C\|u_{0}\|^{2}_{B^{s}_{p,r}}+(\nu^{-1}+T)\int^{T}_{0}\|\psi\|^{2}_{B^{s}_{p,r}(\mathcal{L}^{p})}+\|f\|^{2}_{B^{s-1}_{p,r}}+\|\nabla v\|^{2}_{B^{s}_{p,r}}\|u\|^{2}_{B^{s}_{p,r}}dt+\|\psi\|^{2}_{E^{s}_{p,r}(T)}.
\end{align}
Using Gronwall's inequality, we obtain
\begin{align}
\sup_{t\in[0,T]}\|u\|^{2}_{B^{s}_{p,r}} \leq Ce^{cV(T)}\bigg[\|u_{0}\|^{2}_{B^{s}_{p,r}}+(\nu^{-1}+T)(\varepsilon\|\psi\|^{2}_{E^{s}_{p,r}(T)}+\int^{T}_{0}e^{-cV(t)}(\|\psi\|^{2}_{B^{s}_{p,r}(\mathcal{L}^{p})}+\|f\|^{2}_{B^{s-1}_{p,r}})dt)\bigg].
\end{align}
Now we consider the second esitimate.
Mutiplying both sides of (4.23) by $2^{j(s+1)}$, and taking $l^{r}$-norm with $j>0$, we deduce that
\begin{multline*}
\|(2^{(j+1)s}\|u\|_{L^{p}})_{j\geq 0}\|_{l^{r}}\leq C\bigg(\|(e^{-c\nu t2^{2j}}2^{j}2^{js}\|\Delta_{j}u_{0}\|_{L^{p}})_{j\geq 0}\|_{l^{r}}+\int^{t}_{0}\|(e^{-c\nu (t-t')2^{2j}}2^{j(s+1)}\\(\|div \Delta_{j}\tau\|_{L^{p}}
+\|\Delta_{j}f\|_{L^{p}}+\|\Delta_{j}(v\cdot\nabla u)\|_{L^{p}})_{j>0}\|_{l^{r}}dt')\bigg).
\end{multline*}
Using H\"{o}lder's inequality in $l^{r}$, and by the previous argument for $div~\Delta_{j}\tau$ and $\Delta_{j}(v\cdot\nabla u)$,  we infer that
\begin{multline}
\|(2^{(j+1)s}\|u\|_{L^{p}})_{j\geq 0}\|_{l^{r}}\leq C\bigg(\sup_{j\geq 0}(e^{-c\nu t2^{2j}}2^{j})\|(2^{js}\|\Delta_{j}u_{0}\|_{L^{p}})_{j\geq 0}\|_{l^{r}}+\int^{t}_{0}\sup_{j}(e^{-c\nu (t-t')2^{2j}}2^{2j})\\\|((2^{js}\|\Delta_{j}\tau\|_{L^{p}}
+2^{(j-1)s}\|\Delta_{j}f\|_{L^{p}}+2^{(j-1)s}\|\Delta_{j}(v\otimes u)\|_{L^{p}})_{j>0}\|_{l^{r}}dt')\bigg).
\end{multline}
Note that for any $j\geq 0$, and sufficiently small $\delta>0$
$$e^{-c\nu t2^{2j}}2^{j} \leq \frac{C}{(\nu t)^{\frac{1}{2}-\delta}},~~~~~~e^{-c\nu (t-t')2^{2j}}2^{j} \leq \frac{C}{\nu(t-t')}.$$
Plugging into (4.34), then we have
\begin{multline}
\|(2^{(j+1)s}\|u\|_{L^{p}})_{j\geq 0}\|_{l^{r}}\leq C\bigg(\frac{1}{(\nu t)^{\frac{1}{2}-\delta}}\|(2^{js}\|\Delta_{j}u_{0}\|_{L^{p}})_{j\geq 0}\|_{l^{r}}+\int^{t}_{0}\frac{1}{\nu(t-t')}\\\|((2^{js}\|\Delta_{j}\tau\|_{L^{p}}
+2^{(j-1)s}\|\Delta_{j}f\|_{L^{p}}+2^{js}\|\Delta_{j}(v\otimes u)\|_{L^{p}})_{j>0}\|_{l^{r}}dt')\bigg).
\end{multline}
Taking $L^{2}$-norm over $[0,T]$, and using Young's inequality $\|u*v\|_{L^{2}}\leq C\|u\|_{L^{1}_{w}}\|v\|_{L^{2}},$ we obtain
\begin{multline}
\int^{T}_{0}\|(2^{(j+1)s}\|u\|_{L^{p}})_{j\geq 0}\|^{2}_{l^{r}}dt\leq C\bigg(\nu^{2\delta-1} T^{2\delta} \|(2^{js}\|\Delta_{j}u_{0}\|_{L^{p}})_{j\geq 0}\|_{l^{r}}+\\ \nu^{-1}\int^{T}_{0}\|(2^{js}\|\Delta_{j}\tau\|_{L^{p}}
+2^{(j-1)s}\|\Delta_{j}f\|_{L^{p}}+2^{js}\|\Delta_{j}(v\otimes u)\|_{L^{p}})_{j>0}\|^{2}_{l^{r}}dt')\bigg).
\end{multline}
Taking $L^{2}$-norm over $[0,T]$ for both sides of (4.30), and combinig with the previous argument dealing with $\|\Delta_{j}\tau\|_{L^{p}}$, we deduce that
\begin{multline}
\nu\int^{T}_{0}\|u\|^{2}_{B^{s+1}_{p,r}}dt \leq C (1+\nu T)\bigg(\|u_{0}\|^{2}_{B^{s}_{p,r}}+\varepsilon \|\psi\|^{2}_{E_{T}}+\\ \int^{T}_{0}\|\psi\|^{2}_{B^{s}_{p,r}(\mathcal{L}^{p})}+\|f\|^{2}_{B^{s-1}_{p,r}}+\| v\|^{2}_{B^{s}_{p,r}}\|u\|^{2}_{B^{s}_{p,r}}dt'\bigg).
\end{multline}
Then by the ineqaulity (4.33), we have the following estimate
\begin{align}
\nu\int^{T}_{0}\|u\|^{2}_{B^{s+1}_{p,r}}dt \leq C(1+\nu T)e^{cV(T)}\bigg[\|u_{0}\|^{2}_{B^{s}_{p,r}}+(\varepsilon\|\psi\|^{2}_{E^{s}_{p,r}(T)}+\int^{T}_{0}e^{-cV(t)}(\|\psi\|^{2}_{B^{s}_{p,r}(\mathcal{L}^{p})}+\|f\|^{2}_{B^{s-1}_{p,r}})dt)\bigg].
\end{align}
\end{proof}
\end{lemm}

\begin{rema}
If $u=v$, then Lemma 4.4 holds true for  $V(T)=\int^{T}_{0}\|u\|^{2}_{L^{\infty}}dt$. The proof is similar to that of Lemma 4.4, the only difference is treating with the term $u\cdot \nabla u$. Using the fact that
\begin{align*}
\|u\cdot \nabla u\|_{B^{s-1}_{p,r}}&= \|div(u\otimes u)\|_{B^{s-1}_{p,r}}\leq C\|u\otimes u\|_{B^s_{p,r}}\leq C \|u\|_{L^{\infty}} \|u\|_{B^s_{p,r}},
\end{align*}
then we get the desired result.
\end{rema}

\section{Local well-posedness}
\subsection{Approximate solutions}
~~~~~First, we construct approximate solutions which are smooth (for $x$ variable) solutions of some linear equations.\\ [1ex]
~~~~~Starting for $(u^{0},\psi^{0})\triangleq (S_{0}u_{0},S_{0}\psi_{0})$  we define by induction a sequence $(u^{n},\psi^{n})_{n\in\mathbb{N}}$  by solving the following linear equations:
\begin{align}
\left\{
\begin{array}{ll}
\partial_{t}u^{n+1}+(u^{n}\cdot\nabla)u^{n+1}-\nu\Delta u^{n+1}+\nabla{P^{n+1}}=div~\tau^{n},  div~ u^{n+1}=0,\\[1ex]
\partial_{t}\psi^{n+1}+(u^{n}\cdot\nabla)\psi^{n+1}=div_{R}[-\nabla{u^{n}}R\psi^{n+1}+\psi_{\infty}\nabla_{R}\displaystyle\frac{\psi^{n+1}}{\psi_{\infty}}],  \\[1ex]
\tau^{n+1}_{ij}=\int_{B}R_{i}\otimes\nabla_{j}\mathcal{U}\psi^{n+1} dR, \\[1ex]
u^{n+1}|_{t=0}=S_{n+1}u_{0}, \psi^{n+1}|_{t=0}=S_{n+1}\psi_{0} ,\\[1ex]
\psi_{\infty}\nabla_{R}\displaystyle\frac{\psi^{n+1}}{\psi_{\infty}}=0 ~~~~ on ~~~~ \partial B(0,R_{0}) .\\[1ex]
\end{array}
\right.
\end{align}
Now assume that $u^{n}\in L^{\infty}([0,T];B^{s}_{p,r})\cap L^{2}([0,T];B^{s+1}_{p,r})$ and $\psi^{n} \in C([0,T];B^{s}_{p,r}(\mathcal{L}^{p}))\cap E^{s}_{p,r}(T)$ for some positive $T$. By the previous section's argument we can find a $\psi^{n+1}\in C([0,T];B^{s}_{p,r}(\mathcal{L}^{p}))\cap E^{s}_{p,r}(T)$. Notice that the initial data are smooth, for the linear Stokes equations, there exists a smooth solution $u^{n+1}$. Then Lemma 4.4 guarantees that $u^{n+1}\in L^{\infty}([0,T];B^{s}_{p,r})\cap L^{2}([0,T];B^{s+1}_{p,r})$.

\subsection{Uniform bounds}
Next, we are going to find some positive $T$ such that for which the approximate solutions are uniformly bounded. Define $U^{n}_{1}(t)=\displaystyle\int^{t}_{0}\|u(t')\|^{2}_{B^{s}_{p,r}}dt'$ and $U^{n}_{2}(t)=\displaystyle\int^{t}_{0}\|\nabla u(t')\|^{2}_{B^{s}_{p,r}}dt'$. Then by Lemma 4.2 with $f,g=0$, we have
\begin{align}
\sup_{t\in[0,T]}\|\psi^{n+1}\|_{B^{s}_{p,r}(\mathcal{L}^{p})}+\|\psi^{n+1}\|_{E^{s}_{p,r}(T)}\leq Ce^{T}e^{cU^{n}_{2}(T)}\|\psi_{0}\|_{B^{s}_{p,r}(\mathcal{L}^{p})}.
\end{align}
And by Lemma 4.4 with $f=0$, we obtain
\begin{align}
\sup_{t\in[0,T]}\|u^{n+1}\|^{2}_{B^{s}_{p,r}}\leq Ce^{cU^{n}_{1}(T)}\bigg[\|u_{0}\|^{2}_{B^{s}_{p,r}}+(\nu^{-1}+T)(\varepsilon\|\psi^{n}\|^{2}_{E^{s}_{p,r}(T)}+\int^{T}_{0}e^{-CU^{n}_{1}(t')}\|\psi^{n}\|^{2}_{B^{s}_{p,r}(\mathcal{L}^{p})})dt')\bigg],\\
\nu\int^{T}_{0}\|u^{n+1}\|^{2}_{B^{s+1}_{p,r}}dt \leq C(1+\nu T)e^{cU^{n}_{1}(T)}\bigg[\|u_{0}\|^{2}_{B^{s}_{p,r}}+(\varepsilon\|\psi^{n}\|^{2}_{E^{s}_{p,r}(T)}+\int^{T}_{0}e^{-cU^{n}_{1}(t)}\|\psi^{n}\|^{2}_{B^{s}_{p,r}(\mathcal{L}^{p})}dt)\bigg].
\end{align}
Now fix a $T>0$, such that $$T<\min\bigg\{\nu^{-1},~\ln2,~\frac{1}{C^{2}[\|u_{0}\|^{2}_{B^{s}_{p,r}}+\|\psi_{0}\|^{2}_{B^{s}_{p,r}(\mathcal{L}^{p})}]},~\widetilde{T}\bigg\},$$
where $\widetilde{T}$ denotes the maximal time such that $\frac{4C}{\nu}te^{\frac{2C^{2}}{\nu}A(t)}\leq\frac{1}{2}$, with $$A(t)=\displaystyle\frac{\|u_{0}\|^{2}_{B^{s}_{p,r}}+\|\psi_{0}\|^{2}_{B^{s}_{p,r}(\mathcal{L}^{p})}}{1-C^{2}t(\|u_{0}\|^{2}_{B^{s}_{p,r}}+\|\psi_{0}\|^{2}_{B^{s}_{p,r}(\mathcal{L}^{p})})}.$$
And choose an $\varepsilon$ such that $\frac{4C}{\nu}\varepsilon e^{\frac{2C^{2}}{\nu}A(T)}\leq\frac{1}{2}$.
We claim that for any $n$ and $t\in[0,T]$:
\begin{align}
\|u^{n}(t)\|^{2}_{B^{s}_{p,r}}\leq CA(t), ~~~~\int^{T}_{0}\|u^{n}\|^{2}_{B^{s}_{p,r}}dt\leq \frac{2C}{\nu}A(T), ~~~~\|\psi^{n}\|_{B^{s}_{p,r}(\mathcal{L}^{p})}+\|\psi^{n}\|_{E^{s}_{p,r}(T)}\leq Ce^{\frac{2C^{2}}{\nu}A(T)}\|\psi_{0}\|_{B^{s}_{p,r}}.
\end{align}
By induction, when $n=0$, (5.5) holds true for a fixed $T$. Now suppose that (5.5) is true for $n$.
Plugging (5.5) into (5.2), and using the fact that $\|\nabla u\|_{B^{s}_{p,r}}\leq C\|u\|_{B^{s+1}_{p,r}} $, then we see that
$$\|\psi^{n+1}\|_{B^{s}_{p,r}(\mathcal{L}^{p})}+\|\psi^{n+1}\|_{E^{s}_{p,r}(T)}\leq Ce^{\frac{2C^{2}}{\nu}A(T)}\|\psi_{0}\|_{B^{s}_{p,r}}.$$
Plugging (5.5) into (5.3) we obtain
\begin{align}
\|u^{n+1}(t)\|^{2}_{B^{s}_{p,r}}\leq \frac{C[\|u_{0}\|^{2}_{B^{s}_{p,r}}+(\nu^{-1}+T)(\varepsilon+T)Ce^{\frac{2C^{2}}{\nu}A(T)}\|\psi_{0}\|^{2}_{B^{s}_{p,r}}]}{1-C^{2}t(\|u_{0}\|^{2}_{B^{s}_{p,r}}+\|\psi_{0}\|^{2}_{B^{s}_{p,r}(\mathcal{L}^{p})})}.
\end{align}
The choices of $T$ and $\varepsilon$ ensure that $(\nu^{-1}+T)(\varepsilon+T)Ce^{\frac{2C^{2}}{\nu}A(T)}<1$. Then we deduce that $\|u^{n+1}(t)\|^{2}_{B^{s}_{p,r}}\leq CA(t)$. \\
Plugging (5.5) into (5.4), and by a similar estimate, we have  $\displaystyle\int^{T}_{0}\|u^{n+1}\|^{2}_{B^{s}_{p,r}}dt\leq \frac{2C}{\nu}A(T)$.
Therefore, $u^{n}$ is uniformly bounded in $L^{\infty}([0,T];B^{s}_{p,r})\cap L^{2}([0,T];B^{s+1}_{p,r})$, and $\psi^{n}$ is uniformly bounded in $L^{\infty}([0,T];B^{s}_{p,r}(\mathcal{L}^{p}))\cap E^{s}_{p,r}(T)$.

\subsection{Convergence}
 We are going to show that $(u^{n},\psi^n)$ is a Cauchy sequence in $L^{\infty}([0,T];B^{s-1}_{p,r})\times L^\infty([0,T];B^{s-1}_{p,r}(\mathcal{L}^{p}))\cap E^{s-1}_{p,r}$. By (5.1) we have
 \begin{align}
\left\{
\begin{array}{ll}
\partial_{t}(u^{n+1}-u^{n})+(u^{n}\cdot\nabla)(u^{n+1}-u^{n})-\nu\Delta (u^{n+1}-u^{n})+\nabla(P^{n+1}-P^{n})=div~(\tau^{n}-\tau^{n-1})+F^{n},  \\[1ex]
\partial_{t}(\psi^{n+1}-\psi^{n})+(u^{n}\cdot\nabla)(\psi^{n+1}-\psi^{n})=div_{R}[-\nabla{u^{n}}R(\psi^{n+1}-\psi^{n})+\psi_{\infty}\nabla_{R}\displaystyle\frac{(\psi^{n+1}-\psi^{n})}{\psi_{\infty}}+g^{n}]+f^{n},  \\[1ex]
u^{n+1}-u^{n}|_{t=0}=\Delta_{n}u_{0}, \psi^{n+1}-\psi^{n}|_{t=0}=\Delta_{n}\psi_{0} ,\\[1ex]
\end{array}
\right.
\end{align}
where $F^{n}=-(u^{n}-u^{n-1})\nabla u^{n}$, $g^{n}= -\nabla(u^{n}-u^{n-1})R\psi^{n}$,  $f^{n}=-(u^{n}-u^{n-1})\nabla \psi^{n}$.
By Lemma 4.2 with $\sigma=s-1$, we deuce that
\begin{multline}
\|\psi^{n+1}-\psi^{n}\|^{2}_{B^{s-1}_{p,r}(\mathcal{L}^{p})}+\|\psi^{n+1}-\psi^{n}\|^{2}_{E^{s-1}_{p,r}(T)} \leq Ce^{U^n_2(T)}\bigg(\|\Delta_{n}\psi_{0}\|^{2}_{B^{s-1}_{p,r}(\mathcal{L}^{p})}+\\
\int^T_0\|\nabla(u^{n}-u^{n-1})R\psi^{n}\|^{2}_{B^{s-1}_{p,r}(\mathcal{L}^{p})}
+\|(u^{n}-u^{n-1})\nabla \psi^{n}\|^2_{B^{s-1}_{p,r}(\mathcal{L}^{p})}dt\bigg).
\end{multline}
By a similar calculation as in Lemma 4.4, we obtain
\begin{multline}
\|u^{n+1}-u^{n}\|^{2}_{B^{s-1}_{p,r}}+\int^T_0\|u^{n+1}-u^{n}\|^{2}_{B^{s}_{p,r}}dt \leq C_{T}e^{U^n_{1}(T)}\bigg(\|\Delta_{n}u_{0}\|^{2}_{B^{s-1}_{p,r}}+\varepsilon\|\psi^{n}-\psi^{n-1}\|^{2}_{E^{s-1}_{p,r}(T)}\\
+\int^T_0\|\psi^{n}-\psi^{n-1}\|^{2}_{B^{s-1}_{p,r}(\mathcal{L}^{p})}
+\|(u^{n}-u^{n-1})\nabla u^{n}\|^2_{B^{s-1}_{p,r}}dt\bigg).
\end{multline}
Thanks to the choice of $T$ and that $(u^{n},\psi^{n})$ is uniformly bounded in $L^{\infty}([0,T];B^{s}_{p,r})\cap L^{2}([0,T];B^{s+1}_{p,r})\times L^{\infty}([0,T];B^{s}_{p,r}(\mathcal{L}^{p}))\cap E^{s}_{p,r}(T)$, we get a constant $C_T$ independent of $n$, such that
\begin{align}
\|u^{n+1}-u^{n}\|^{2}_{L^\infty_T(B^{s-1}_{p,r})}+\int^T_0\|u^{n+1}-u^{n}\|^{2}_{B^{s}_{p,r}}dt &\leq C_{T}\bigg(\|\Delta_{n}u_{0}\|^{2}_{B^{s-1}_{p,r}}+\varepsilon\|\psi^{n}-\psi^{n-1}\|^{2}_{E^{s-1}_{p,r}(T)}\\
&+\int^T_0\|\psi^{n}-\psi^{n-1}\|^{2}_{B^{s-1}_{p,r}(\mathcal{L}^{p})}
\nonumber+\|\nabla u^{n}\|^{2}_{B^{s-1}_{p,r}}\|(u^{n}-u^{n-1})\|^2_{B^{s-1}_{p,r}}dt\bigg)\\
\nonumber&\leq \nonumber C_{T}\bigg(\|\Delta_{n}u_{0}\|^{2}_{B^{s-1}_{p,r}}+\varepsilon\|\psi^{n}-\psi^{n-1}\|^{2}_{E^{s-1}_{p,r}(T)}\\
\nonumber &+\int^T_0\|\psi^{n}-\psi^{n-1}\|^{2}_{B^{s-1}_{p,r}(\mathcal{L}^{p})}
+\|(u^{n}-u^{n-1})\|^2_{B^{s-1}_{p,r}}dt\bigg).
\end{align}
And
\begin{multline}
\|\psi^{n+1}-\psi^{n}\|^{2}_{L^\infty_T(B^{s-1}_{p,r}(\mathcal{L}^{p}))}+\|\psi^{n+1}-\psi^{n}\|^{2}_{E^{s-1}_{p,r}(T)} \leq C_T\bigg(\|\Delta_{n}\psi_{0}\|^{2}_{B^{s-1}_{p,r}(\mathcal{L}^{p})}\\
 +\int^T_0\|\nabla(u^{n}-u^{n-1})\|^2_{B^{s-1}_{p,r}}\|\psi^{n}\|^{2}_{B^{s-1}_{p,r}(\mathcal{L}^{p})}
+\|(u^{n}-u^{n-1})\|^2_{B^{s-1}_{p,r}}\|\nabla \psi^{n}\|^2_{B^{s-1}_{p,r}(\mathcal{L}^{p})}dt\bigg) \\
 \leq C_T\bigg(\|\Delta_{n}\psi_{0}\|^{2}_{B^{s-1}_{p,r}(\mathcal{L}^{p})}
+\int^T_0\|(u^{n}-u^{n-1})\|^2_{B^{s}_{p,r}}
+\|(u^{n}-u^{n-1})\|^2_{B^{s-1}_{p,r}}dt\bigg).
\end{multline}
By a direct calculation, we obtain
\begin{align}
\|\Delta_{n}u_0\|^{2}_{B^{s-1}_{p,r}}&=(\sum_{|j-n|\leq1}2^{jr(s-1)}\|\Delta_j\Delta_{n}u_0\|^r_{L^p})^\frac{2}{r}\\
\nonumber &\leq C 2^{2ns}\|\Delta_{n}u_0\|^{2}_{L^p} 2^{-2n}\leq C 2^{-2n} \|u_0\|^2_{B^s_{p,r}}.
\end{align}
By a similar argument, we have $\|\Delta_{n}\psi_0\|^{2}_{B^{s}_{p,r}}\leq C2^{-2n}\|\psi_0\|^{2}_{B^{s}_{p,r}}$. If we define
\begin{align}
A_n(T)=\|u^{n+1}-u^{n}\|^{2}_{L^\infty_T(B^{s-1}_{p,r})}+\int^T_0\|u^{n+1}-u^{n}\|^{2}_{B^{s}_{p,r}}dt,  \\
B_n(T)=\|\psi^{n+1}-\psi^{n}\|^{2}_{L^\infty_T(B^{s-1}_{p,r}(\mathcal{L}^{p}))}+\|\psi^{n+1}-\psi^{n}\|^{2}_{E^{s-1}_{p,r}(T)}.
\end{align}
Plugging (5.12) into (5.10), we have
\begin{align}
A_n(T)\leq C_T [2^{-2n}\|u_0\|^2_{B^s_{p,r}}+(\varepsilon+T)(A_{n-1}(T)+B_{n-1}(T))], \\
B_n(T)\leq C_T (2^{-2n}\|\psi_0\|^2_{B^s_{p,r}(\mathcal{L}^{p})}+A_{n-1}(T) ).
\end{align}
Plugging (5.16) into (5.15), we deduce that
\begin{align}
A_n(T)\leq C_T 2^{-2n}(\|u_0\|^2_{B^s_{p,r}}+\|\psi_0\|^2_{B^s_{p,r}(\mathcal{L}^{p})})+C^2_T(\varepsilon+T)(A_{n-1}(T)+A_{n-2}(T)).
\end{align}
Now let $\varepsilon$ and $T$ be sufficiently small such that $C^2_T(\varepsilon+T)<\frac{1}{4}$, then
\begin{align}
A_n(T)\leq C_T 2^{-2n}(\|u_0\|^2_{B^s_{p,r}}+\|\psi_0\|^2_{B^s_{p,r}(\mathcal{L}^{p})})+\frac{1}{4}(A_{n-1}(T)+A_{n-2}(T)).
\end{align}
If $n<0$, we may set $A_n=0$, (5.18) still holds true. Then we have
\begin{align}
\sum^m_{n=1}A_n(T)\leq C_T \sum^m_{n=1}2^{-2n}(\|u_0\|^2_{B^s_{p,r}}+\|\psi_0\|^2_{B^s_{p,r}(\mathcal{L}^{p})})+\frac{1}{4}\sum^m_{n=1}(A_{n-1}(T)+A_{n-2}(T)).
\end{align}
Since $A_n(T)>0$ and $A_n=0$ if $n<0$, it follows that
$$\sum^m_{n=1}(A_{n-1}(T)+A_{n-2}(T))\leq \frac{1}{2}\sum^m_{n=1}A_n(T).$$
Then
\begin{align}
\sum^m_{n=1}A_n(T)\leq C_T \sum^m_{n=1}2^{-2n}(\|u_0\|^2_{B^s_{p,r}}+\|\psi_0\|^2_{B^s_{p,r}(\mathcal{L}^{p})}).
\end{align}
This implies that $\sum^\infty_{n=1}A_n(T)$ is convergent. Hence, by (5.16), $\sum^\infty_{n=1}B_n(T)$ is also convergent. Then we deduce that $(u^{n},\psi^n)$ is a Cauchy sequence in $L^{\infty}([0,T];B^{s-1}_{p,r})\times L^\infty([0,T];B^{s-1}_{p,r}(\mathcal{L}^{p}))\cap E^{s-1}_{p,r}$. Thus, there exists $(u,\psi)\in L^{\infty}([0,T];B^{s-1}_{p,r})\times L^\infty([0,T];B^{s-1}_{p,r}(\mathcal{L}^{p}))$ such that
$$ u^n\rightarrow u~~~in~~L^{\infty}([0,T];B^{s-1}_{p,r})~~and~~\psi^n\rightarrow \psi~~~in~~L^\infty([0,T];B^{s-1}_{p,r}(\mathcal{L}^{p})).$$
Since $u^n$ and $\psi^n$ are uniform bounded in $L^{\infty}([0,T];B^{s}_{p,r})\cap L^2([0,T];B^{s+1}_{p,r})\times L^\infty([0,T];B^{s}_{p,r}(\mathcal{L}^{p}))$. The Fatou property for Besov spaces ensures that $(u,\psi)\in L^{\infty}([0,T];B^{s}_{p,r})\cap L^2([0,T];B^{s+1}_{p,r})\times L^\infty([0,T];B^{s}_{p,r}(\mathcal{L}^{p}))$. An interpolation argument ensures that the convergence holds true for any $s'<s$. Passing to the limit in (5.1) in the weak sense, we conclude that  $(u,\psi)$ is indeed a solution of (1.2).
\subsection{Regularity}
 ~~~Now we check that $(u,\psi)\in C([0,T];B^{s}_{p,r})\times C([0,T];B^{s}_{p,r}(\mathcal{L}^{p}))$, when $r$ is finite, and $(u,\psi)\in C_w([0,T];B^{s}_{p,\infty})\times C_w([0,T];B^{s}_{p,\infty}(\mathcal{L}^{p})).$
 Proposition 4.3  guarantees that $\psi$ is in the desired space.\\
 ~~~Since $u\in L^\infty([0,T];B^s_{p,r})$, it follows that $u\cdot\nabla u\in L^\infty([0,T];B^{s-1}_{p,r})$ and  $\Delta u \in L^\infty([0,T];B^{s-2}_{p,r})$. Combining with $\psi\in L^\infty([0,T];B^{s}_{p,r}(\mathcal{L}^{p}))\cap E^s_{p,r}(T)$, we deduce that $div~\tau\in L^\infty([0,T];B^{s-1}_{p,r}).$ Applying $div$ to (1.2) we obtain
 $$ \Delta P=div[~(u\nabla u)+div~\tau],$$
 from which we deduce that $$\nabla P=\nabla (\Delta)^{-1}div~[~(u\nabla u)+div~\tau].$$
 If $p<\infty$,  we have
 $$\|\nabla P\|_{B^{s-1}_{p,r}}\leq C \|u\nabla u\|_{B^{s-1}_{p,r}}+\|div~\tau\|_{B^{s-1}_{p,r}}.$$
 So from the equations (1.1) we obtain $\partial_t u\in L^\infty([0,T];B^{s-2}_{p,r}),$ then $u\in C([0,T];B^{s-2}_{p,r}).$ An interpolation argument ensures that $u\in C([0,T];B^{s'}_{p,r})$ for any $s'<s$.
 Notice that $u\in L^\infty([0,T];B^{s}_{p,r})$. If $r<\infty$, for any $t_1,t_2$, for any $\varepsilon>0$, there exists $N$ such that
 \begin{align}
 \sup_{t\in[0,T]}(\sum_{j\geq N}2^{jsr}\|u(t)\|^r_{L^p})^\frac{1}{r} \leq \frac{\varepsilon}{4}.
 \end{align}
 Hence
 \begin{align*}
 \|u(t_1)-u(t_2)\|_{B^{s}_{p,r}}&\leq  (\sum_{-1\leq j< N}2^{jsr}\|u(t_1)-u(t_2)\|^r_{L^p})^\frac{1}{r}+2\sup_{t\in[0,T]}(\sum_{j\geq N}2^{jsr}\|u(t)\|^r_{L^p})^\frac{1}{r} \\
 &\leq  (\sum_{-1\leq j< N}2^{jsr}\|u(t_1)-u(t_2)\|^r_{L^p})^\frac{1}{r} + \frac{\varepsilon}{4} \\
 &\leq  2^N(\sum_{-1\leq j< N}2^{j(s-1)r}\|u(t_1)-u(t_2)\|^r_{L^p})^\frac{1}{r}+\frac{\varepsilon}{4}\\
 &\leq  2^N\|u(t_1)-u(t_2)\|_{B^s_{p,r}}+\frac{\varepsilon}{4},
 \end{align*}
 from which we deduce that $u\in C([0,T];B^s_{p,r})$.
 If $r=\infty$, for any $\phi\in S$
 $$\langle u(t),\phi\rangle=\langle S_j u(t),\phi\rangle+\langle (Id-S_j) u(t),\phi\rangle=\langle S_j u(t),\phi\rangle+\langle u(t),(Id-S_j) \phi\rangle.$$
 Since $S_j u(t)\in B^{s'}_{p,r}$, it follows that $\langle S_j u(t),\phi\rangle\in C[0,T]$ and $\|(Id-S_j)\phi\|_{L^\infty}\rightarrow 0,~~j\rightarrow\infty$. Then  $\langle S_j u(t),\phi\rangle$ uniformly converges to $\langle u(t),\phi\rangle$. So $\langle u(t),\phi\rangle$ is continuous, which implies that $u\in C_w([0,T];B^s_{p,\infty})$.

\subsection{Uniqueness}
Assume that $(u,\psi)$ and $(v,\phi)$ are two solutions of (1.2) with the same initial data. Then we have
 \begin{align}
\left\{
\begin{array}{ll}
\partial_{t}(u-v)+(v\cdot\nabla)(u-v)-\nu\Delta (u-v)+\nabla(P_1-P_2)=div~(\tau_1-\tau_2)+F,  \\[1ex]
\partial_{t}(\psi-\phi)+(v\cdot\nabla)(\psi-\phi)=div_{R}[-\nabla{u^{n}}R(\psi-\phi)+\psi_{\infty}\nabla_{R}\displaystyle\frac{(\psi-\phi)}{\psi_{\infty}}+g]+f,  \\[1ex]
u-v|_{t=0}=0, \psi-\phi|_{t=0}0 ,\\[1ex]
\end{array}
\right.
\end{align}
where $F=-(u-v)\nabla v$, $g= -\nabla(u-v)R\psi$,  $f=-(u-v)\nabla \phi$, and $P_1$ corresponds to $u$, $\tau_1$ corresponds to $\psi$, $P_2$ corresponds to $v$, $\tau_2$ corresponds to $\phi$ respectively. By a similar calculation as in Section 5.3, we have
\begin{align}
\|u-v\|^{2}_{L^\infty_T(B^{s-1}_{p,r})}+\int^T_0\|u-v\|^{2}_{B^{s}_{p,r}}dt \leq C_{T}\bigg(\varepsilon\|\psi-\phi\|^{2}_{E^{s-1}_{p,r}(T)}
+\int^T_0\|\psi-\phi\|^{2}_{B^{s-1}_{p,r}(\mathcal{L}^{p})}
+\|u-v\|^2_{B^{s-1}_{p,r}}dt\bigg),
\end{align}
\begin{align}
\|\psi-\phi\|^{2}_{L^\infty_T(B^{s-1}_{p,r}(\mathcal{L}^{p}))}+\|\psi-\phi\|^{2}_{E^{s-1}_{p,r}(T)}
 \leq C_T\bigg(\int^T_0\|(u-v)\|^2_{B^{s}_{p,r}}+\|(u-v)\|^2_{B^{s-1}_{p,r}}dt\bigg).
\end{align}
Plugging (5.24) into (5.23), we deduce that
\begin{align}
\|u-v\|^{2}_{L^\infty_T(B^{s-1}_{p,r})}+\int^T_0\|u-v\|^{2}_{B^{s}_{p,r}}dt \leq C^2_{T}(\varepsilon+T)\bigg(
\|u-v\|^2_{L^\infty_T(B^{s-1}_{p,r})}+\int^T_0\|u-v\|^{2}_{B^{s}_{p,r}}dt\bigg).
\end{align}
So if $\varepsilon$ and $T$ are small enough such that $C^2_{T}(\varepsilon+T)<1$, we obtain $u=v~~for~~a.e.~~(t,x)$.
By the inequality (5.24), we have $\psi=\phi ~~for~~a.e.~~(t,x,R)$.
Splitting the interval $[0,T]$ into $[0,T_1]$, $[T_1,T_2]$,...,$[T_k,T_{k+1}]$, where each subinterval satisfies $C^2_{T}(\varepsilon+T)<1$. Then we see that $(u,\psi)=(v,\phi)$ in the interval $[0,T_1]$. Since $u(T_1)=v(T_1),~~\psi(T_1)=\phi(T_1)$, we deduce that $(u,\psi)=(v,\phi)$ in the interval $[T_1,T_2]$. Repeating the argument we have proved the uniqueness.

\section{Blow-up criterion}
In this section we give the proof of the blow-up criterion for (1.2).\\
\textbf{Proof of Theorem 2.2}:
If $T^{*}<\infty$, we have
$$\|\psi\|^2_{L^\infty_{T^{*}}(B^{s}_{p,r}(\mathcal{L}^{p}))}+\|\psi\|^2_{E^s_{p,r}(T*)}=\infty.$$
Since $\|\psi\|^2_{L^\infty_{T^{*}}(B^{s}_{p,r}(\mathcal{L}^{p}))}+\|\psi\|^2_{E^s_{p,r}(T*)}<\infty$, it follows that $T^{*}$ is not the maximal time.
Now assume that $$\int^{T^*}_0\|u\|^2_{L^\infty}<\infty.$$
Using Remark 4.5 with $f=0$, we deduce that
\begin{align}
\sup_{t\in[0,T^{*})}\|u\|^{2}_{B^{s}_{p,r}}&\leq M(T^*)   \\
&\triangleq Ce^{c\int^{T^*}_0\|u\|^2_{L^\infty}dt}\bigg[\|u_{0}\|^{2}_{B^{s}_{p,r}}+
(\varepsilon\|\psi\|^{2}_{E^{s}_{p,r}(T^*)}+\int^{T^*}_{0}\|\psi\|^{2}_{B^{s}_{p,r}(\mathcal{L}^{p})})dt')\bigg].
\end{align}
By the assumption, we have that $M(T^*)<\infty$. Let $\delta$ be small enough such that
$$\delta<\min\{\nu^{-1}, \ln2,\frac{1}{C^2[M(T^*)+\|\psi\|^2_{L^\infty_{T^{*}}{(B^{s}_{p,r}(\mathcal{L}^{p}))}}]}
\}.$$
Then by the argument as in Section 5.1, we have a solution $\widetilde{u}$ of (1.1) with initial data $u(T^*-\frac{\delta}{2})$. By the uniqueness, we deduce that $\widetilde{u}(t) = u(t + T^{*}-\frac{\delta}{2})$ on $[0,\frac{\delta}{2})$.
So the solution $\widetilde{u}$ extends the solution $u$ beyond $T^*$. This contradicts the fact that $T^*$ is the lifespan.

\section{Global existence for small data}
~~~~~   In this section, we proved that the solution is global in time if the initial data is close to equilibrium $(0,\psi_\infty)$. Firstly we need the following Poincar$\acute{e}$ inequality with weight. The proof is similar as in Proposition 3.4 in \cite{Masmoudi.W}.
   \begin{lemm}
   If $\widetilde{\psi}$ satisfy $\int_B \widetilde{\psi}dR=0$ and
   $\displaystyle\int_B\bigg|\nabla _R  \bigg(\displaystyle\frac{\widetilde{\psi} }{\psi _\infty }\bigg)^{\frac{p}{2}}\bigg|^2{\psi _\infty} dR<\infty$ with $p\geq 2$. There exists a constant $C$ such that
   $$\int_B \bigg|\frac{\widetilde{\psi}}{\psi_\infty}\bigg|^p \psi_\infty dR\leq C \displaystyle\int_B\bigg|\nabla _R  \bigg(\displaystyle\frac{\widetilde{\psi} }{\psi _\infty }\bigg)^{\frac{p}{2}}\bigg|^2{\psi _\infty} dR.$$
   \begin{proof}
   We argue by contradiction. Assume that the sequence $\widetilde{\psi}_n$ satisfies:
   \begin{align}
   \int_B \widetilde{\psi}_ndR=0,~~~\int_B \bigg|\frac{\widetilde{\psi}_n}{\psi_\infty}\bigg|^p \psi_\infty dR=1,~~~\int_B\bigg|\nabla _R  \bigg(\displaystyle\frac{\widetilde{\psi}_n }{\psi _\infty }\bigg)^{\frac{p}{2}}\bigg|^2{\psi _\infty} dR\rightarrow 0.
   \end{align}
   Hence, $\sqrt{\psi_\infty}\nabla _R  \bigg(\displaystyle\frac{\widetilde{\psi}_n }{\psi _\infty }\bigg)^{\frac{p}{2}}$ tends to $0$ in $L^2(B)$ and $\nabla _R  \bigg(\displaystyle\frac{\widetilde{\psi}_n }{\psi _\infty }\bigg)^{\frac{p}{2}}$ tends to $0$ in $L^2_{loc}(B)$. Then we deduce that $  \bigg(\displaystyle\frac{\widetilde{\psi}_n }{\psi _\infty }\bigg)^{\frac{p}{2}}$ tends to some constant $c$ in $L^2_{loc}(B)$. Thus we obtain $\displaystyle\frac{\widetilde{\psi}_n }{\psi _\infty }$ tends to some constant $c$ in $L^p_{loc}(B)$. Since $\widetilde{\psi}_n$ is bounded in $L^p(B)$, it follows that $\int_B \widetilde{\psi}_n dR$ tends to some constant $c$. Using the fact that $\int_B \widetilde{\psi}_ndR=0$, we infer that $c=0$. There exists a subsequence $\widetilde{\psi}_{n_k}$ such that $\widetilde{\psi}_{n_k}\rightarrow 0 $ almost every.
   Denote that $x=1-|R|$. By a similar calculation as in Section 3, we have $\psi_\infty \sim  x^k$. By the Hardy-type inequality (for more details, one can refer to Section 3.2 in \cite{Masmoudi.G}), we deduce that for some $\beta>0$
   $$\int_B \bigg|\frac{\widetilde{\psi}_n}{x^\beta\psi_\infty}\bigg|^p \psi_\infty dR<C\bigg[\int_B\bigg|\displaystyle\frac{\widetilde{\psi}_n }{\psi _\infty }\bigg|^p\psi_\infty dR+\int_B \bigg|\nabla _R  \bigg(\displaystyle\frac{\widetilde{\psi}_n }{\psi _\infty }\bigg)^{\frac{p}{2}}\bigg|^2\bigg]<C.$$
   This gives some tightness of the sequence of $\bigg|\displaystyle\frac{\widetilde{\psi}_n }{\psi _\infty }\bigg|^p$, thus we have
   $$\lim_{n\rightarrow\infty}\int_B\bigg|\displaystyle\frac{\widetilde{\psi}_n }{\psi _\infty }\bigg|^p\psi_\infty dR=\int_B\lim_{n\rightarrow\infty}\bigg|\displaystyle\frac{\widetilde{\psi}_n }{\psi _\infty }\bigg|^p\psi_\infty dR=0,$$
   which contradicts (7.1).
   \end{proof}
   \end{lemm}

   Now we are going to prove Theorem 2.3. Denote that $\widetilde{\psi}=\psi-\psi_\infty $. Since $(0,\psi_\infty)$ is the equilibrium of (1.2), it follows that $(u,\widetilde{\psi})$ is also a solution of (1.2).\\
   \textbf{Proof of Theorem 2.3:} We argue by contradiction, assume that the lifespan $T^{*}$ is finite. Firstly we claim that there exists a constant $M$ independent with $T$, such that
   \begin{multline}
   \sup_{t\in[0,T^{*})}\|u(t)\|_{B^{s}_{p,r}}+\nu\|u\|_{L^2_{T^{*}}(B^{s+1}_{p,r})}+\sup_{t\in[0,T^{*})}\|\widetilde{\psi}(t)\|_{B^{s}_{p,r}(\mathcal{L}^{p})}+\|\widetilde{\psi}\|_{E^{s}_{p,r}(T^{*})} \\ \leq M(\|u_0\|_{B^{s}_{p,r}}+\|\psi_0-\psi_\infty\|_{B^{s}_{p,r}(\mathcal{L}^{p})}).
   \end{multline}
   If $T$ is small enough, by the argument as in local well-posedness, we have
   \begin{multline}
   \sup_{t\in[0,T]}\|u(t)\|_{B^{s}_{p,r}}+\nu\|u\|_{L^2_{T}(B^{s+1}_{p,r})}+\sup_{t\in[0,T]}\|\widetilde{\psi}(t)\|_{B^{s}_{p,r}(\mathcal{L}^{p})}+\|\widetilde{\psi}\|_{E^{s}_{p,r}(T)} \\ \leq M(\|u_0\|_{B^{s}_{p,r}}+\|\psi_0-\psi_\infty\|_{B^{s}_{p,r}(\mathcal{L}^{p})}).
   \end{multline}
   Thus we can define $\overline{T}$  as follow:
   \begin{multline}
   \overline{T}=\sup\bigg\{T:\sup_{t\in[0,T]}\|u(t)\|_{B^{s}_{p,r}}+\nu\|u\|_{L^2_{T}(B^{s+1}_{p,r})}+\sup_{t\in[0,T]}\|\widetilde{\psi}(t)\|_{B^{s}_{p,r}(\mathcal{L}^{p})}+\|\widetilde{\psi}\|_{E^{s}_{p,r}(T)} \\ \leq M(\|u_0\|_{B^{s}_{p,r}}+\|\psi_0-\psi_\infty\|_{B^{s}_{p,r}(\mathcal{L}^{p})}) \bigg\}.
   \end{multline}
   We seek to prove that  $\overline{T}=T^{*}$, then the claim holds true. If $\overline{T}<T^{*}$, for any $t<\overline{T}$. By a similar calculation as in Lemma 4.2 with $f,g=0$, we deduce that

   \begin{multline}
\frac{1}{p}\partial_{t}\int_{\mathbb{R}^{d}\times B}\bigg|\frac{\Delta_{j}\widetilde{\psi}}{\psi_{\infty}}\bigg|^{p}\psi_{\infty}dxdR+\frac{2(p-1)}{p^{2}}\int_{\mathbb{R}^{d}\times B}\bigg|\nabla_{R}\bigg(\frac{\Delta_{j}\widetilde{\psi}}{\psi_{\infty}}\bigg)^{\frac{p}{2}}\bigg|^{2}\psi_{\infty}dxdR \leq \\ C\int_{\mathbb{R}^{d}\times B}\bigg(\Delta_{j}\bigg(\nabla{u}R\frac{\widetilde{\psi}}{\psi_{\infty}}\bigg)\bigg)^{2}\bigg|\frac{\Delta_{j}\psi}{\psi_{\infty}}\bigg|^{p-2}\psi_{\infty}dxdR
+\int_{\mathbb{R}^{d}\times B} R^1_{j}\bigg|\frac{\Delta_{j}\widetilde{\psi}}{\psi_{\infty}}\bigg|^{p-1}dxdR,
\end{multline}
where $R^1_j=[u\cdot\nabla,\Delta_j]\psi$.  \\
Using H\"{o}lder's inequality, and by the definition of Besov spaces we have
   \begin{multline}
\partial_{t}\|\Delta_j\widetilde{\psi}\|^p_{L^{p}_{x}(\mathcal{L}^{p})}+\frac{2(p-1)}{p}\int_{\mathbb{R}^{d}\times B}\bigg|\nabla_{R}\bigg(\frac{\Delta_{j}\widetilde{\psi}}{\psi_{\infty}}\bigg)^{\frac{p}{2}}\bigg|^{2}\psi_{\infty}dxdR \leq \\ pC\bigg(c^2_j2^{-2js}\|\nabla u R \widetilde{\psi}\|^2_{ B^{s}_{p,r}(\mathcal{L}^{p})}\|\Delta_j\widetilde{\psi}\|^{p-2}_{L^{p}_{x}(\mathcal{L}^{p})}
+\|R^1_j\|_{L^{p}_{x}(\mathcal{L}^{p})}\|\Delta_j\widetilde{\psi}\|^{p-1}_{L^{p}_{x}(\mathcal{L}^{p})}\bigg).
\end{multline}
By Lemmas 3.15 and 3.17, we obtain
   \begin{multline}
\partial_{t}\|\Delta_j\widetilde{\psi}\|^p_{L^{p}_{x}(\mathcal{L}^{p})}+\frac{2(p-1)}{p}\int_{\mathbb{R}^{d}\times B}\bigg|\nabla_{R}\bigg(\frac{\Delta_{j}\widetilde{\psi}}{\psi_{\infty}}\bigg)^{\frac{p}{2}}\bigg|^{2}\psi_{\infty}dxdR \leq \\ pC\bigg(c^2_j2^{-2js}\|\nabla u\|^2_{ B^{s}_{p,r}} \| \widetilde{\psi}\|^2_{ B^{s}_{p,r}(\mathcal{L}^{p})}\|\Delta_j\widetilde{\psi}\|^{p-2}_{L^{p}_{x}(\mathcal{L}^{p})}
+c_j2^{-js}\|u\|_{ B^{s}_{p,r}}\|\widetilde{\psi}\|_{{ B^{s}_{p,r}}(\mathcal{L}^{p})}\|\Delta_j\widetilde{\psi}\|^{p-1}_{L^{p}_{x}(\mathcal{L}^{p})}\bigg).
\end{multline}
Let $C_p$ denote a constant dependent on $p$. By Lemma 7.1, we deduce that
   \begin{multline}
\partial_{t}\|\Delta_j\widetilde{\psi}\|^p_{L^{p}_{x}(\mathcal{L}^{p})}+C_p\|\Delta_j\widetilde{\psi}\|^p_{L^{p}_{x}(\mathcal{L}^{p})} \leq \\ pC\bigg(c^2_j2^{-2js}\|\nabla u\|^2_{ B^{s}_{p,r}} \| \widetilde{\psi}\|^2_{ B^{s}_{p,r}(\mathcal{L}^{p})}\|\Delta_j\widetilde{\psi}\|^{p-2}_{L^{p}_{x}(\mathcal{L}^{p})}
+c_j2^{-js}\|u\|_{ B^{s}_{p,r}}\|\widetilde{\psi}\|_{{ B^{s}_{p,r}}(\mathcal{L}^{p})}\|\Delta_j\widetilde{\psi}\|^{p-1}_{L^{p}_{x}(\mathcal{L}^{p})}\bigg).
\end{multline}
Thus we obtain
 \begin{multline}
\partial_{t}\|\Delta_j\widetilde{\psi}\|^r_{L^{p}_{x}(\mathcal{L}^{p})}+C_p\|\Delta_j\widetilde{\psi}\|^r_{L^{p}_{x}(\mathcal{L}^{p})} \leq \\ pC\bigg(c^2_j2^{-2js}\|\nabla u\|^2_{\dot B^{s}_{p,r}} \| \widetilde{\psi}\|^2_{ B^{s}_{p,r}(\mathcal{L}^{p})}\|\Delta_j\widetilde{\psi}\|^{r-2}_{L^{p}_{x}(\mathcal{L}^{p})}
+c_j2^{-js}\|u\|_{ B^{s}_{p,r}}\|\widetilde{\psi}\|_{{ B^{s}_{p,r}}(\mathcal{L}^{p})}\|\Delta_j\widetilde{\psi}\|^{r-1}_{L^{p}_{x}(\mathcal{L}^{p})}\bigg).
\end{multline}
Multiplying both sides of (7.9) by $2^{jrs}$, taking sum of $j$ form $-1$ to $\infty$, and using  H\"{o}lder's inequality, we deduce that
   \begin{align}
\partial_{t}\|\widetilde{\psi}\|^r_{B^{s}_{p,r}(\mathcal{L}^{p})}+C_p\|\widetilde{\psi}\|^r_{ B^{s}_{p,r}(\mathcal{L}^{p})} \leq  pC\bigg(\|\nabla u\|^2_{B^{s}_{p,r}} \|\widetilde{\psi}\|^{r}_{ B^{s}_{p,r}(\mathcal{L}^{p})}
+\|u\|_{B^{s}_{p,r}}\|\widetilde{\psi}\|^{r}_{B^{s}_{p,r}(\mathcal{L}^{p})}\bigg).
\end{align}
Since $t<\overline{T}$, it follows that
$$\|u\|_{B^s_{p,r}}\leq M\big(\|u_0\|_{B^s_{p,r}}+\|\psi_0-\psi_\infty\|_{B^s_{p,r}(\mathcal{L}^{p})}\big)\leq Mc_0.$$
So if $p C M c_0<\frac{C_p}{2}$, we get
 \begin{align}
\partial_{t}\|\widetilde{\psi}\|^r_{ B^{s}_{p,r}(\mathcal{L}^{p})}+\frac{C_p}{2}\|\widetilde{\psi}\|^r_{ B^{s}_{p,r}(\mathcal{L}^{p})} \leq  pC\|\nabla u\|^2_{ B^{s}_{p,r}} \|\widetilde{\psi}\|^{r}_{ B^{s}_{p,r}(\mathcal{L}^{p})},
\end{align}
or
 \begin{align}
\partial_{t}\|\widetilde{\psi}\|_{ B^{s}_{p,r}(\mathcal{L}^{p})}+\frac{C_p}{2}\|\widetilde{\psi}\|_{ B^{s}_{p,r}(\mathcal{L}^{p})} \leq  pC\|\nabla u\|^2_{ B^{s}_{p,r}} \|\widetilde{\psi}\|_{ B^{s}_{p,r}(\mathcal{L}^{p})}.
\end{align}
Integrating over $[0,t]$ with respect to $t$, we have
  \begin{align*}
\|\widetilde{\psi}(t)\|_{ B^{s}_{p,r}(\mathcal{L}^{p})} &\leq \|\psi_0-\psi_\infty\|_{ B^{s}_{p,r}(\mathcal{L}^{p})}+ pC\int^{\overline{T}}_0\|\nabla u\|^2_{ B^{s}_{p,r}} \|\widetilde{\psi}\|_{ B^{s}_{p,r}(\mathcal{L}^{p})}dt \\
& \leq  \|\psi_0-\psi_\infty\|_{ B^{s}_{p,r}(\mathcal{L}^{p})}+ pCMc_0 \sup_{t\in[0,\overline{T}]}\|\widetilde{\psi}\|_{ B^{s}_{p,r}(\mathcal{L}^{p})}.
\end{align*}
Taking $L^{\infty}$-norm for both sides of the above inequality, if $p C M c_0<\frac{1}{2}$, we obtain
\begin{align}
\sup_{t\in[0,\overline{T}]}\|\widetilde{\psi}(t)\|_{ B^{s}_{p,r}(\mathcal{L}^{p})} &\leq 2\|\psi_0-\psi_\infty\|_{ B^{s}_{p,r}(\mathcal{L}^{p})}.
\end{align}
Integrating (7.7) over $[0,\overline{T}]$ with respect to $t$, and multiplying $2^{jps}$ for both sides of (7.7) and then taking the $l^{\frac{r}{p}}$-norm, we deduce that
\begin{align*}
\|\widetilde{\psi}\|^p_{E^s_{p,r}(\overline{T})}&\leq \frac{2p}{p-1}\|\psi_0-\psi_\infty\|^p_{B^{s}_{p,r}(\mathcal{L}^{p})}+C_p(\int^{\overline{T}}_{0}\|\nabla u\|^2_{B^s_{p,r}}\|\widetilde{\psi}\|^p_{B^{s}_{p,r}(\mathcal{L}^{p})}dt+\int^{\overline{T}}_{0}\| u\|_{B^s_{p,r}}\|\widetilde{\psi}\|^p_{B^{s}_{p,r}(\mathcal{L}^{p})}dt)  \\
&\leq \frac{2p}{p-1}\|\psi_0-\psi_\infty\|^p_{B^{s}_{p,r}(\mathcal{L}^{p})}+C_p (M c_0)^p \sup_{t\in[0,\overline{T}]}\|\widetilde{\psi}(t)\|^p_{ B^{s}_{p,r}(\mathcal{L}^{p})} + C_p M c_0 \int^{\overline{T}}_{0}\|\widetilde{\psi}\|^p_{B^{s}_{p,r}(\mathcal{L}^{p})}dt.
\end{align*}
By Lemma 7.1, we deduce that $\displaystyle\int^{\overline{T}}_{0}\|\widetilde{\psi}\|^p_{B^{s}_{p,r}(\mathcal{L}^{p})}dt\leq C\|\widetilde{\psi}\|^p_{E^s_{p,r}(\overline{T})}$. If $C C_p M c_0\leq \frac{1}{2} $, we infer that
\begin{align}
\|\widetilde{\psi}\|^p_{E^s_{p,r}(\overline{T})}&\leq \frac{4p}{p-1}\|\psi_0-\psi_\infty\|^p_{B^{s}_{p,r}(\mathcal{L}^{p})}+2C_p(M c_0)^p \sup_{t\in[0,\overline{T}]}\|\widetilde{\psi}(t)\|^p_{ B^{s}_{p,r}(\mathcal{L}^{p})}.
\end{align}
So plugging (7.13) into (7.14) and choosing $c_0$ small enough we have
\begin{align}
\|\widetilde{\psi}\|^p_{E^s_{p,r}(\overline{T})}&\leq \frac{8p}{p-1}\|\psi_0-\psi_\infty\|^p_{B^{s}_{p,r}(\mathcal{L}^{p})}.
\end{align}
Since $p\in [2,\infty)$, it follows that
\begin{align}
\|\widetilde{\psi}\|_{E^s_{p,r}(\overline{T})}\leq (\frac{8p}{p-1})^{\frac{1}{p}}\|\psi_0-\psi_\infty\|_{B^{s}_{p,r}(\mathcal{L}^{p})}\leq 4\|\psi_0-\psi_\infty\|_{B^{s}_{p,r}(\mathcal{L}^{p})}.
\end{align}
For the Navier-Stokes equations, we can write that
\begin{align}
u=e^{t\nu\Delta}u_0+\int^t_0 e^{(t-t')\nu\Delta}(div\widetilde{\tau}-u\nabla u-\nabla P) dt'.
\end{align}
 Using the fact that $s>0$, $B^{s}_{p,r}=\dot B^{s}_{p,r} \cap L^{p}$. Applying $\dot\Delta_j$ to (7.17) we obtain
\begin{align}
\dot\Delta_j u=e^{t\nu\Delta}\dot\Delta_j u_0+\int^t_0 e^{(t-t')\nu\Delta}(div\dot\Delta_j\widetilde{\tau}-\dot\Delta_j(u\nabla u)-\nabla \dot\Delta_j P) dt'.
\end{align}
By Lemma 3.2, we deduce that
\begin{align}
\|\dot\Delta_j u\|_{L^p}\leq Ce^{-tc\nu2^{2j}}\|\dot\Delta_j u_0\|_{L^p}+\int^t_0 e^{-(t-t')c\nu2^{2j}}(\|div\dot\Delta_j\widetilde{\tau}\|_{L^p}+\|\dot\Delta_j(u\nabla u)\|_{L^p}+\|\nabla \dot\Delta_j P)\|_{L^p})dt'.
\end{align}
By a similar argument as in Lemma 4.1, we deduce that $\|\nabla \dot\Delta_j P)\|_{L^p}\leq C(\|div\dot\Delta_j\widetilde{\tau}\|_{L^p}+\|\dot\Delta_j(u\nabla u)\|_{L^p})$. Thus
\begin{align}
\|\dot\Delta_j u\|_{L^p}\leq e^{-tc\nu2^{2j}}\|\dot\Delta_j u_0\|_{L^p}+C\int^t_0 e^{-(t-t')c\nu2^{2j}}(\|div\dot\Delta_j\widetilde{\tau}\|_{L^p}+\|\dot\Delta_j(u\nabla u)\|_{L^p})dt'.
\end{align}
Using Young's inequality, we get
\begin{align*}
\|\dot\Delta_j u\|_{L^p} &\leq \|\dot\Delta_j u_0\|_{L^p}+C(\int^t_0 2^{-2j}(\|div\dot\Delta_j\widetilde{\tau}\|^2_{L^p}+\|\dot\Delta_j(u\nabla u)\|^2_{L^p})dt')^{\frac{1}{2}} \\
&\leq \|\dot\Delta_j u_0\|_{L^p}+ C(\int^t_0 (\|\dot\Delta_j\widetilde{\tau}\|^2_{L^p}+2^{-2j}\|\dot\Delta_j(u\nabla u)\|^2_{L^p})dt')^{\frac{1}{2}}.
\end{align*}
By the definition of Besov spaces, we have
\begin{align}
\|\dot\Delta_j u\|^2_{L^p} \leq \|\dot\Delta_j u_0\|^2_{L^p}+C(\int^t_0 (\|\dot\Delta_j\widetilde{\tau}\|^2_{L^p}+c^2_j2^{-2js}\|(u\nabla u)\|^2_{\dot B^{s-1}_{p,r}})dt'),
\end{align}
 where $c_j$ denotes an element of the unit sphere of $l^r$. Multiplying both sides of (7.21) by $2^{2js}$, taking the $l^{\frac{r}{2}}$-norm, and by Lemma 7.1 and Corollary 3.13 we deduce that
\begin{align}
\|u\|^2_{\dot B^s_{p,r}} \leq \|u_0\|^2_{\dot B^s_{p,r}}+C(\|\widetilde{\psi}\|^2_{\dot E^s_{p,r}(T)}+\int^{t}_0\|(u\nabla u)\|^2_{\dot B^s_{p,r}}dt').
\end{align}
Since $E^{s}_{p,r}(T)\hookrightarrow \dot E^s_{p,r}(T)$ and  $B^s_{p,r}\hookrightarrow\dot B^s_{p,r}$, it follows that
\begin{align*}
\|u\|^2_{\dot B^s_{p,r}} &\leq \|u_0\|^2_{\dot B^s_{p,r}}+16\|\psi_0-\psi_\infty\|^2_{B^{s}_{p,r}(\mathcal{L}^{q})}+\int^{t}_0\|(\nabla u)\|^2_{B^s_{p,r}}dt'\sup_{t\in\overline{T}}\|u\|^2_{\dot B^s_{p,r}} ) \\
& \leq \|u_0\|^2_{\dot B^s_{p,r}}+16\|\psi_0-\psi_\infty\|^2_{B^{s}_{p,r}(\mathcal{L}^{q})}+CMc_0\sup_{t\in\overline{T}}\|u\|^2_{\dot B^s_{p,r}} ).
\end{align*}
If $CMc_0<\frac{1}{2}$, then we obtain
\begin{align}
\sup_{t\in\overline{T}}\|u\|^2_{\dot B^s_{p,r}}\leq 2\|u_0\|^2_{\dot B^s_{p,r}}+32\|\psi_0-\psi_\infty\|^2_{B^{s}_{p,r}(\mathcal{L}^{q})}.
\end{align}
Hence
\begin{align}
\sup_{t\in\overline{T}}\|u\|_{\dot B^s_{p,r}}\leq 2\|u_0\|_{\dot B^s_{p,r}}+8\|\psi_0-\psi_\infty\|_{B^{s}_{p,r}(\mathcal{L}^{q})}.
\end{align}
Multiplying both sides of (7.19) by $2^{2j(s+1)}$ and taking the $L^2([0,\overline{T}])$-norm, by a similar calculation, we deduce that
\begin{align}
\nu\|u\|_{L^2_{\overline{T}}(\dot B^s_{p,r})}\leq 2\|u_0\|_{\dot B^s_{p,r}}+8\|\psi_0-\psi_\infty\|_{B^{s}_{p,r}(\mathcal{L}^{q})}.
\end{align}
Now we estimate the $L^p$-norm. By (1.2), we write that
\begin{align}
\partial_{t}u^k+(u\cdot\nabla)u^k-\nu\Delta{u^k}+\partial_k{P}=div\widetilde{\tau^k},
\end{align}
where $u^k$ is the $k$ component of $u$. Note that $div u=0$. Multiplying $sgn (u^k)|u^k|^{p-1}$ both sides of (7.26) and integrating over $\mathbb{R}^d$, we have
\begin{align}
\frac{1}{p}\partial_{t}\int_{\mathbb{R}^d}|u^k|^pdx+\nu(p-1)(\int_{\mathbb{R}^d}|\nabla{u^k}|^2|u^k|^{p-2}dx=(p-1)\int_{\mathbb{R}^d}{P}|\partial_k u^k||u^k|^{p-2}dx-\int_{\mathbb{R}^d}\widetilde{\tau^k}\nabla{u^k}|u^k|^{p-2}dx).
\end{align}
By  Cauchy-Schwarz's inequality, we deduce that
 \begin{align}
\frac{1}{p}\partial_{t}\int_{\mathbb{R}^d}|u^k|^pdx+\frac{\nu(p-1)}{2}\int_{\mathbb{R}^d}|\nabla{u^k}|^2|u^k|^{p-2}dx\leq C_p(\int_{\mathbb{R}^d}P^2|u^k|^{p-2}dx+\int_{\mathbb{R}^d}|\widetilde{\tau^k}|^2|u^k|^{p-2}dx).
\end{align}
Using H\"{o}lder's inequality, we obtain
 \begin{align}
\partial_{t}\|u\|^p_{L^p}+\frac{\nu(p-1)}{2}\int_{\mathbb{R}^d}|\nabla{u}|^2|u|^{p-2}dx\leq C_p(\|P\|^2_{L^p}\|u\|^{p-2}_{L^p}+\|\widetilde{\tau}\|^2_{L^p}\|u\|^{p-2}_{L^p}).
\end{align}
Since $P=\Delta^{-1}div (u\nabla u+div \widetilde{\tau})=\Delta^{-1}div div (u\otimes u+\widetilde{\tau})$, it follows that
$$\|P\|_{L^p}\leq C(\|u\otimes u\|_{L^p}+\|\widetilde{\tau}\|_{L^p})\leq C\|u\|_{L^\infty}\|u\|_{L^p}+\|\widetilde{\tau}\|_{L^p}.$$
By Lemma 7.1 and Corollary 3.13, we get
\begin{align}
\|\widetilde{\tau}\|^p_{L^p}\leq  C\int_{\mathbb{R}^d\times B}\bigg|\nabla_{R}\bigg(\displaystyle\frac{\widetilde{\psi}}{\psi_\infty}\bigg)^{\frac{p}{2}}\bigg|^2 \psi_{\infty} dxdR \leq C\sum_{j\geq-1}\int_{\mathbb{R}^d\times B}\bigg|\nabla_{R}\bigg(\displaystyle\frac{\widetilde{\Delta_{j}\psi}}{\psi_\infty}\bigg)^{\frac{p}{2}}\bigg|^2 \psi_{\infty} dxdR.
\end{align}
Plugging into (7.29), we deduce that
\begin{align}
\nonumber \partial_{t}\|u\|^p_{L^p}&\leq C_p(\|u\|^2_{L^\infty}\|u\|^{p}_{L^p}+2\|\widetilde{\tau}\|^2_{L^p}\|u\|^{p-2}_{L^p}) \\
&\leq C_p(\|u\|^2_{L^\infty}\|u\|^{p}_{L^p} + 2C\sum_{j\geq-1}\int_{\mathbb{R}^d\times B}\bigg|\nabla_{R}\bigg(\displaystyle\frac{\widetilde{\Delta_{j}\psi}}{\psi_\infty}\bigg)^{\frac{p}{2}}\bigg|^2 \psi_{\infty} \|u\|^{p-2}_{L^p}) \\
&\leq C_p(\|u\|^2_{L^\infty}\|u\|^{p}_{L^p} + 2C\sum_{j\geq-1}2^{-js}2^{js}\int_{\mathbb{R}^d\times B}\bigg|\nabla_{R}\bigg(\displaystyle\frac{\widetilde{\Delta_{j}\psi}}{\psi_\infty}\bigg)^{\frac{p}{2}}\bigg|^2 \psi_{\infty} \|u\|^{p-2}_{L^p}).
\end{align}
Note that $s>0$. Integrating over $[0,t]$ with respect to $t$ and using H\"{o}lder's inequality, we obtain
\begin{align}
\|u\|^p_{L^p}&\leq \|u_0\|^p_{L^p}+ C_p(\int^{\overline{T}}_0\|u\|^2_{L^\infty}dt+\|\widetilde{\psi}\|_{E^s_{p,r}(T)}) \sup_{t\in[0,\overline{T}]}\|u\|^{p}_{L^p}.
\end{align}
Since $B^s_{p,r}\hookrightarrow L^\infty$, it follows that
 \begin{align}
\|u\|^p_{L^p}&\leq \|u_0\|^p_{L^p}+ C_p Mc_0 \sup_{t\in[0,\overline{T}]}\|u\|^{p}_{L^p}.
\end{align}
If $C_p Mc_0<\frac{1}{2}$, then we have
\begin{align}
\sup_{t\in[0,\overline{T}]}\|u\|^p_{L^p}&\leq 2\|u_0\|^p_{L^p} ~~or~~\sup_{t\in[0,\overline{T}]}\|u\|_{L^p}\leq 2\|u_0\|_{L^p}.
\end{align}
Combining with (7.13), (7.16), (7.24), (7.25), (7.35) and using the fact that $B^s_{p,r}=\dot B^s_{p,r}\cap L^p$, we get
\begin{multline}
  \sup_{t\in[0,\overline{T}]}\|u(t)\|_{B^{s}_{p,r}}+\nu\|u\|_{L^2_{\overline{T}}(B^{s+1}_{p,r})}+\sup_{t\in[0,\overline{T}]}\|\widetilde{\psi}(t)\|_{B^{s}_{p,r}(\mathcal{L}^{p})}+\|\widetilde{\psi}\|_{E^{s}_{p,r}(\overline{T})} \\ \leq 16(\|u_0\|_{B^{s}_{p,r}}+\|\psi_0-\psi_\infty\|_{B^{s}_{p,r}(\mathcal{L}^{p})}).
   \end{multline}
   Note that $u\in C([0,T^*);B^{s}_{p,r}) ,\widetilde{\psi} \in C([0,T^*);B^{s}_{p,r}(\mathcal{L}^{p})) $. If we set $M=32$, then we can find a $T_1\in[\overline{T},T^{*}]$, such that
\begin{multline}
  \sup_{t\in[0,T_1]}\|u(t)\|_{B^{s}_{p,r}}+\nu\|u\|_{L^2_{T_1}(B^{s+1}_{p,r})}+\sup_{t\in[0,T_1]}\|\widetilde{\psi}(t)\|_{B^{s}_{p,r}(\mathcal{L}^{p})}+\|\widetilde{\psi}\|_{E^{s}_{p,r}(\overline{T})} \\ \leq 32(\|u_0\|_{B^{s}_{p,r}}+\|\psi_0-\psi_\infty\|_{B^{s}_{p,r}(\mathcal{L}^{p})})=M(\|u_0\|_{B^{s}_{p,r}}+\|\psi_0-\psi_\infty\|_{B^{s}_{p,r}(\mathcal{L}^{p})}).
   \end{multline}
   This contradicts the definition of $\overline{T}$. Thus we have $\overline{T}=T^{*}$. Therefore we obtain
   \begin{multline}
    \sup_{t\in[0,T^{*})}\|u(t)\|_{B^{s}_{p,r}}+\nu\|u\|_{L^2_{T^{*}}(B^{s+1}_{p,r})}+\sup_{t\in[0,T^{*})}\|\widetilde{\psi}(t)\|_{B^{s}_{p,r}(\mathcal{L}^{p})}+\|\widetilde{\psi}\|_{E^{s}_{p,r}(T^{*})} \\ \leq M(\|u_0\|_{B^{s}_{p,r}}+\|\psi_0-\psi_\infty\|_{B^{s}_{p,r}(\mathcal{L}^{p})}).
   \end{multline}
   Now set $T_{\delta}=T^{*}-\frac{\delta}{2}$. Then we can construct a solution $(\widehat{u},\widehat{\psi})$ with initial data $\|u(T_{\delta})\|_{B^{s}_{p,r}}$ and $\|\widetilde{\psi}\|_{E^{s}_{p,r}(T_{\delta})}$.
   This implies the solution can be extended outside $[0,T^{*})$, which contradicts the lifespan $T^{*}$. Thus the solution is global. This completes the proof.\\

\smallskip
\noindent\textbf{Acknowledgments} This work was partially supported by
NNSFC (No. 11271382 and No. 10971235), RFDP (No. 20120171110014),
and the key project of Sun Yat-sen University. 

\phantomsection
\addcontentsline{toc}{section}{\refname}

\end{document}